\documentclass[a4paper, BCOR=0.0mm, DIV=calc]{scrartcl}
\usepackage[T1]{fontenc}
\usepackage[utf8]{inputenc}
\usepackage[ngerman]{babel}
\usepackage[osf,sc]{mathpazo}
\usepackage{textcomp}
\usepackage{microtype}
\usepackage{graphicx, color}
\usepackage{amsmath, amssymb, amsfonts}
\usepackage{booktabs}
\usepackage{nicefrac}
\usepackage{wasysym,pifont}
\usepackage{tikz}
\usepackage{graphicx}
\usepackage{float}
\usepackage{wrapfig}
\usetikzlibrary{calc,intersections}
\KOMAoptions{twoside=false, twocolumn=false, headinclude=false, footinclude=false, mpinclude=false, pagesize=auto}
\recalctypearea

\setkomafont{section}{\mdseries\scshape\Large}
\addtokomafont{subsection}{\mdseries\scshape}
\setkomafont{title}{\mdseries\scshape}

\begin{document}
\title{Observations on continued fractions
\footnote{
Original title: "' De fractionibus continuis Observationes"', first published in "`\textit{Commentarii academiae scientiarum Petropolitanae 11}, 1750, pp. 32-81"', reprinted in  "`\textit{Opera Omnia}: Series 1, Volume 14, pp. 291 - 349 "', Eneström-Number E123, translated by: Alexander Aycock, for the "`Euler-Kreis Mainz"'}}
\author{Leonhard Euler}
\date{ }
\maketitle

\paragraph*{§1}

After  I had started to investigate continued fractions in the last year and to develop this almost completely new branch of analysis, in the mean time I made several observations, which will be helpful for expanding this theory even further. Therefore, because the exploration of this theory seems to be very useful for all branches of Analysis in general, I will take on this subject again, and explain everything in more detail. So let this continued fraction be propounded

\[
A+\cfrac{B}{C+\cfrac{D}{E+\cfrac{F}{G+\cfrac{H}{I+\text{etc.}}}}}
\]
whose true value will be found by continuing the following  series to infinity

\[
A+\frac{B}{1P}-\frac{BD}{PQ}+\frac{BDF}{QR}-\frac{BDFH}{RS}-\text{etc.},
\]
in which series the letters $P$, $Q$, $R$, $S$ etc. obtain the following values

\[
P=C, \quad Q=EP+D, \quad R=GQ+FP, \quad S=IR+HQ \quad \text{etc.}
\] 
But this series  always converges, no matter how the letters $B$, $C$, $D$, $E$, $F$ etc. either grow or decrease, as long they are all positive; for, any arbitrary term is smaller than the preceding one, but greater than the following; this is immediately seen by the rule used to form the values $P$, $Q$, $R$, $S$ etc.

\paragraph*{§2}

If vice versa this infinite series was propounded

\[
\frac{B}{P}-\frac{BD}{PQ}+\frac{BDF}{QR}-\frac{BDFH}{RS}+\text{etc.},
\]
its sum can conveniently be expressed by the following continued fraction. Because it is

\[
C=P, \quad E=\frac{Q-D}{P}, \quad G=\frac{R-FP}{Q}, \quad I=\frac{S-HQ}{R} \quad \text{etc.},
\]
one will have this continued fraction equal to that series

\[
\cfrac{B}{P+\cfrac{D}{\cfrac{Q-D}{P}+\cfrac{F}{\cfrac{R-FP}{Q}+\cfrac{H}{\cfrac{S-HQ}{R}+\cfrac{K}{\text{etc.}}}}}}
\]
or

\[
\cfrac{B}{P+\cfrac{DP}{Q-D+\cfrac{FPQ}{R-FP+\cfrac{HQR}{S-HQ+\cfrac{KRS}{\text{etc.}}}}}}
\]
If this series was given

\[
\frac{a}{p}-\frac{b}{q}+\frac{c}{r}-\frac{d}{s}+\frac{e}{t}-\text{etc.},
\]
because of

\[
B=a, \quad D=b:a, \quad F=c:b, \quad H=d:c, \quad K=e:d \quad \text{etc.}
\]
and

\[
P=p, \quad Q=q:p, \quad R=pr:q, \quad S=qs:pr, \quad T=prt:qs \quad \text{etc.}
\]
the sum of this series

\[
\frac{a}{p}-\frac{b}{q}+\frac{c}{r}-\frac{d}{s}+\frac{e}{t}-\text{etc.}
\]
will be equal to the following continued fraction

\[
\cfrac{a}{p+\cfrac{b:a}{\cfrac{aq-bp}{app}+\cfrac{c:b}{\cfrac{p^2(br-cq)}{bqq}+\cfrac{d:c}{\cfrac{q^2(cs-dr)}{cp^2r^2}+\cfrac{e:d}{\cfrac{p^2r^2(dt-es)}{dq^2s^2}+\text{etc.}}}}}}
\]
\[
=\cfrac{a}{p+\cfrac{bp^2}{aq-bp+\cfrac{acqq}{br-cq+\cfrac{bdrr}{cs-dr+\cfrac{cess}{dt-es+\text{etc.}}}}}}
\]

\paragraph*{§3}

To illustrate these things by some examples, let us take this series

\[
1-\frac{1}{2}+\frac{1}{3}-\frac{1}{4}+\frac{1}{5}-\frac{1}{6}+\text{etc.},
\]
whose sum is $=\log 2$ or $\int \frac{dx}{1+x}$, if after the integration $x$ is put $=1$; therefore, it will be be

\[
a=b=c=d=\text{etc.}=1, \quad p=1, \quad q=2, \quad r=3, \quad s=4 \quad \text{etc.}
\]
and

\[
p=1, \quad aq-bp=1, \quad br-cq=1, \quad cs-dr=1 \quad \text{etc.}
\]
Hence it will be

\[
\int \frac{dx}{1+x}=\cfrac{1}{1+\cfrac{1}{1+\cfrac{4}{1+\cfrac{9}{1+\cfrac{16}{1+\text{etc.}}}}}}
\]
or the value of this continued fraction is $=\log 2$.

\paragraph*{§4}

Now let us consider this series 

\[
1-\frac{1}{3}+\frac{1}{5}-\frac{1}{7}+\frac{1}{9}-\text{etc.},
\]
whose sum is the area of the circle, whose diameter is $=1$ or $=\int \frac{dx}{1+x^2}$ having put $x=1$ after the integration. Therefore, it will be

\[
a=b=c=d=\text{etc.}=1 \quad \text{and} \quad p=1, \quad q=3, \quad r=5, \quad s=7 \quad \text{etc.},
\]
whence it is

\[
\int \frac{dx}{1+xx}=\cfrac{1}{1+\cfrac{1}{2+\cfrac{9}{2+\cfrac{25}{2+\cfrac{49}{2+\text{etc.}}}}}}
\]
which is \textsc{Brouncker's} continued fraction he exhibited for the quadrature of the circle.

\paragraph*{§5}

In like manner by taking other series of this kind, the following conversions of integral formulas into continued fractions will result, after having put $x=1$ after the integration:

\begin{alignat*}{9}
&\int \frac{dx}{1+x^3}&&=\cfrac{1}{1+\cfrac{1^2}{3+\cfrac{4^2}{3+\cfrac{7^2}{3+\cfrac{10^2}{3+\text{etc.}}}}}} \quad \quad && &\int \frac{dx}{1+x^4}&&=\cfrac{1}{1+\cfrac{1^2}{4+\cfrac{5^2}{4+\cfrac{9^2}{4+\cfrac{13^2}{4+\text{etc.}}}}}}\\
&\int \frac{dx}{1+x^5}&&=\cfrac{1}{1+\cfrac{1^2}{5+\cfrac{6^2}{5+\cfrac{11^2}{5+\cfrac{16^2}{5+\text{etc.}}}}}} \quad \quad && &\int \frac{dx}{1+x^6}&&=\cfrac{1}{1+\cfrac{1^2}{6+\cfrac{7^2}{6+\cfrac{13^2}{6+\cfrac{19^2}{6+\text{etc.}}}}}}
\end{alignat*}

\paragraph*{§6}

Hence it follows that it will be in general

\[
\int \frac{dx}{1+x^m}=\cfrac{1}{1+\cfrac{1^2}{m+\cfrac{(m+1)^2}{m+\cfrac{(2m+1)^2}{m+\cfrac{(3m+1)^2}{m+\text{etc.}}}}}}
\]
having put $x=1$ after the integration. And if $m$ was a fractional number, one will have

\[
\int \frac{dx}{1+x^{\frac{m}{n}}}=\cfrac{1}{1+\cfrac{n}{m+\cfrac{(m+n)^2}{m+\cfrac{(2m+n)^2}{m+\cfrac{(3m+n)^2}{m+\text{etc.}}}}}}
\]

\paragraph*{§7}

Now let us consider the formula $\int \frac{x^{n-1}dx}{1+x^m}$, which integrated and having put $x=1$ afterwards, gives this series

\[
\frac{1}{n}-\frac{1}{m+n}+\frac{1}{2m+n}-\frac{1}{3m+n}+\text{etc.}
\]
Hence it will be

\[
a=b=c=d=\text{etc.}=1 \quad \text{and} \quad p=n, \quad q=m+n, \quad r=2m+n, \quad s=3m+n \quad \text{etc.}
\]
Hence one will have 

\[
\int \frac{x^{n-1}dx}{1+x^m}=\cfrac{1}{n+\cfrac{n^2}{m+\cfrac{(m+n)^2}{m+\cfrac{(2m+n)^2}{m+\text{etc.}}}}}
\]
which continued fraction coincides with the last one found.

\paragraph*{§8}

Now let this formula be propounded $\int \frac{x^{n-1}dx}{(1+x^m)^{\frac{\mu}{\nu}}}$, which integrated and having put $x=1$ afterwards gives this series

\[
\frac{1}{n}-\frac{\mu}{\nu(m+n)}+\frac{\mu(\mu +\nu)}{1 \cdot 2 \nu^2(2m+n)}-\frac{\mu(\mu +\nu)(\mu+ 2 \nu)}{1 \cdot 2 \cdot 3 \nu^3 (3m+n)}+\text{etc.},
\]
which compared to the general one yields

\[
a=1, \quad b=\mu, \quad c=\mu(\mu+\nu), \quad d=\mu(\mu+\nu)(\mu+2 \nu) \quad \text{etc.},
\]
\[
p=n, \quad q=\nu (m+n), \quad R=2\nu^2(2m+n), \quad s=6 \nu^3(3m+n),
\]
\[
t=24\nu^4(4m+n) \quad \text{etc.}
\]
and

\begin{alignat*}{9}
&aq&&-bp&&=\nu m+(\nu-\mu)n,\\
&br&&-cq&&=\mu \nu (3 \nu-\mu)m+\mu \nu(\nu -\mu)n,\\
&cs&&-dr&&=2 \mu \nu^3(\mu +\nu)(m(5\nu -2 \mu)+n(\nu- \mu)),\\
&dt&&-es&&=6 \mu \nu^3(\mu+\nu)(\mu+2\nu)(m(7 \nu-3 \mu)+n(\nu-\mu))\\
&  &&   && \quad \quad \quad \quad \quad \quad \quad \quad \text{etc.};
\end{alignat*}
having substituted these values and after some simplifications one will have

\[
\int \frac{x^{n-1}dx}{(1+x^m)^{\frac{\mu}{\nu}}}
\]
\begin{small}
\[
= \cfrac{1}{n+\cfrac{\mu n^2}{\nu m+(\nu-\mu)n+\cfrac{\nu(\mu+\nu)(m+n)^2}{(3\nu-\mu)+(\nu-\mu)n+\cfrac{2\nu(\mu+2 \nu)(2m+n)^2}{(5\nu-2\mu)m+(\nu-\mu)n+\cfrac{3\nu(\mu+3\nu)(3m+n)^2}{(7\nu-3\mu)m+(\nu-\mu)n+\text{etc.}}}}}}
\]
\end{small}
Let $\mu=1$ and $\nu=2$; it will be

\[
\int \frac{x^{n-1}dx}{\sqrt{1+x^m}}=\cfrac{1}{n+\cfrac{n^2}{2m+n+\cfrac{6(m+n)^2}{5m+n+\cfrac{20(2m+n)^2}{8m+n+\cfrac{42(3m+n)^2}{11m+n+\cfrac{72(4m+n)^2}{14m+n+\text{etc.}}}}}}}
\]

\paragraph*{§9}

But if $\nu=1$ and $\mu$ was an integer, the following continued fractions will result

\[
\int \frac{x^{n-1}dx}{(1+x^m)^2}=\cfrac{1}{n+\cfrac{2n^2}{m-n+\cfrac{1 \cdot 3 (m+n)^2}{m-n+\cfrac{2 \cdot 4(2m+n)^2}{m-n+\cfrac{3 \cdot 5 (3m+n)^2}{m-n+\cfrac{4 \cdot 6(4m+n)^2}{m-n+\text{etc.}}}}}}}
\]
\[
\int \frac{x^{n-1}dx}{(1+x^m)^3}= \cfrac{1}{n+\cfrac{3n^3}{m-2n+\cfrac{1 \cdot 4(m+n)^2}{-2n+\cfrac{2 \cdot 5 (2m+n)^2}{-m-2n+\cfrac{3 \cdot 6(3m+n)^2}{-2m-2n+\cfrac{4 \cdot 7(4m+n)^2}{-3m-2n+\text{etc.}}}}}}}
\]
which expressions,  as the following ones, because of the negative quantities, do not converge, but diverge.

\paragraph*{§10}

All these things follow from the conversion of the general continued fraction given in § $1$ into an infinite series

\[
A+\frac{B}{1P}-\frac{BD}{PQ}+\frac{BDF}{QR}-\frac{BDFH}{RS}+\text{etc.}
\]
But this same series, by adding each two terms, is transformed into this one

\[
A+\frac{BE}{1Q}+\frac{BDFI}{QS}+\frac{BDFHKN}{SV}+\text{etc.}
\]
But it is

\[
C=P=\frac{Q-D}{E}, \quad G=\frac{S-HQ}{IQ}-\frac{F(Q-D)}{EQ}, \quad L=\frac{V-MS}{NS}-\frac{K(S-HQ)}{IS} \quad \text{etc.}
\]
Therefore, this infinite series

\[
A+\frac{BE}{Q}+\frac{BDFI}{QS}+\frac{BDFHKN}{SV}+\text{etc.}
\]
is converted into the following continued fraction

\[
A+\cfrac{B}{\cfrac{Q-D}{E}+\cfrac{D}{E+\cfrac{F}{\cfrac{E(S-HQ)-FI(Q-D)}{EIQ}+\cfrac{H}{I+\cfrac{K}{\cfrac{I(V-MS)-KN(S-HQ)}{INS}+\text{etc.}}}}}};
\]
and having cleared all the fractions this expression becomes

\[
A+\cfrac{BE}{Q-D+\cfrac{D}{1+\cfrac{FIQ}{E(S-HQ)-FI(Q-D)+\cfrac{EHQ}{1+\cfrac{KNS}{I(V-MS)-KN(S-HQ)+\cfrac{IMS}{1+\text{etc.}}}}}}}
\]

\paragraph*{§11}

If now vice versa this infinite series is propounded
\[
\frac{a}{p}+\frac{b}{q}+\frac{c}{r}+\frac{d}{s}+\frac{e}{t}+\text{etc.}
\]
and compared  to the preceding, it will be

\[
Q=p, \quad S=\frac{q}{p},\quad V=\frac{pr}{q}, \quad X=\frac{qs}{pr}, \quad Z=\frac{prt}{qs} \quad \text{etc.}
\]
and similarly

\[
E=\frac{a}{B},\quad I=\frac{b}{BDF},\quad N=\frac{c}{BDFHK} \quad \text{etc.},
\]
by means of which values the propounded series is converted into this continued fraction

\[
\cfrac{a}{p-D+\cfrac{D}{1+\cfrac{bp:1}{Da\bigg(\cfrac{q}{p}-Hp\bigg)-b(p-D)+\cfrac{DHap:1}{1+\cfrac{cq:p}{Hb\bigg(\cfrac{pr}{q}-\cfrac{Mq}{p}\bigg)-c\bigg(\cfrac{q}{p}-Hp\bigg)+\cfrac{HMbq:p}{1+\cfrac{dp:q}{Mc\bigg(\frac{qs}{pr}-\text{etc.}}}}}}}};
\]
this continued fraction contains innumerable new quantities, which were not contained in the propounded series.

\paragraph*{§12}

But because from § $2$ this series 

\[
\frac{b}{p}-\frac{bd}{pq}+\frac{bdf}{qr}-\frac{bdfh}{rs}+\text{etc.}
\]
is equal to this continued fraction

\[
\cfrac{b}{p+\cfrac{dp}{q-d+\cfrac{fpq}{r-fp+\cfrac{hqr}{s-hq+\cfrac{krs}{\text{etc.}}}}}},
\]
if this last series is reduced to the preceding, it will be

\[
b=BE, \quad d=\frac{-DFI}{E}, \quad f=\frac{-HKN}{I} \quad \text{etc.},
\]
\[
p=Q, \quad q=S, \quad r=V, \quad s=X \quad \text{etc.}
\]
Hence the continued fraction given in the preceding paragraph is transformed into this one

\[
A+\cfrac{BE}{Q-\cfrac{DFI \cdot Q}{ES+DFI-\cfrac{EHKN \cdot QS}{IV+HKNQ-\cfrac{IMOR \cdot SV}{NX+MORS+\text{etc.}}}}}
\]
whose structure is easily understood.

\paragraph*{§13}

But that series

\[
A+\frac{B}{P}-\frac{BD}{PQ}+\frac{BDF}{QR}-\text{etc.},
\]
 we  found  from the general continued fraction first is easily transformed into this form

\[
A+\frac{B}{2P}+\frac{BE}{2Q}-\frac{BDG}{2PR}+\frac{BDFI}{2QS}-\frac{BDFHL}{2RT}+\text{etc.};
\]
this form, if the letters $C$, $E$, $G$, $I$ etc. are expressed by means of the given equations, becomes 

\[
A+\frac{B}{2P}+\frac{B(Q-D)}{2PQ}-\frac{BD(R-FP)}{2PQR}+\frac{BDF(S-HQ)}{2QRS}-\text{etc.};
\]
hence the following continued fraction is equal to the series
\[
A+\cfrac{B}{P+\cfrac{DP}{Q-D+\cfrac{FPQ}{R-FP+\cfrac{HQR}{A-HQ+\text{etc.}}}}}
\]

\paragraph*{§14}

Therefore, all these things  follow  directly from the consideration of continued fractions and I  have mentioned the most observations of this kind already in the preceding dissertation\footnote{Euler refers to his paper "De fractionibus continuis dissertatio". This is E71 in the Eneström-Index}. But now having treated these things I will consider others and will explain several ways so to get to continued fractions as to assign the values of continued fractions of that kind by means of integrations. Therefore at first, because the  expression given by \textsc{Brouncker} for the quadrature of the circle is not only demonstrated, but was also found a priori, I want to investigate other similar expressions, either found by \textsc{Brouncker} himself or by \textsc{Wallis}; for, they were listed up by \textsc{Wallis}, but it was not indicated clearly enough, whether \textsc{Brouncker} found all of them or merely that one  exhibited for the quadrature of the circle. But after this,  I will also prove those remaining continued fractions, which seem to be  more profound, using completely different principles and will teach how to  find a lot more continued fractions of this kind.

\paragraph*{§15}

But everything concerning continued fractions in \textsc{Wallis's} book reduces to this that the product of the following two continued fractions is $=a^2$:

\[
a-1+\cfrac{1}{2(a-1)+\cfrac{9}{2(a-1)+\cfrac{25}{2(a-1)+\text{etc.}}}}
\]
and

\[
a+1+\cfrac{1}{2(a+1)+\cfrac{9}{2(a+1)+\cfrac{25}{2(a+1)+\text{etc.}}}}
\]
But because likewise it is

\[
(a+2)^2=
a+1+\cfrac{1}{2(a+1)+\cfrac{9}{2(a+1)+\text{etc.}}} \quad \times \quad  a+3+\cfrac{1}{2(a+3)+\cfrac{9}{2(a+3)+\text{etc.}}}
\]
by proceeding in this way to infinity one will find

\[
a \cdot \frac{~~~~~~~~~a(a+4)}{(a+2)(a+2)}\frac{(a+8)(a+~~8)}{(a+6)(a+10)}\frac{(a+12)(a+12)}{(a+10)(a+14)}~~\text{etc.}
\]
\[
=a-1+\cfrac{1}{2(a-1)+\cfrac{9}{2(a-1)+\cfrac{25}{2(a-1)+\text{etc.}}}}
\]

\paragraph*{§16}

If this product consisting of an infinite number of factors is examined by the method given in the preceding paper\footnote{This time Euler refers to his paper "De productis ex infinitis factoribus ortis". This is E122 in the Eneström-Index.}, one will find 

\[
\frac{a(a+4)(a+4)(a+8)\text{etc.}}{(a+2)(a+2)(a+6)(a+6)\text{etc.}}=\frac{\int x^{a+1}dx:\sqrt{1-x^4}}{\int x^{a-1}dx:\sqrt{1-x^4}}.
\]
Therefore, the value of this continued fraction

\[
=a-1+\cfrac{1}{2(a-1)+\cfrac{9}{2(a-1)+\cfrac{25}{2(a-1)+\text{etc.}}}}
\]
will become equal to this expression

\[
a \frac{\int x^{a+1}dx:\sqrt{1-x^4}}{\int x^{a-1}dx:\sqrt{1-x^4}}
\]
having put $x=1$ after each integration.

\paragraph*{§17}

This theorem by which the value of a  far-extending continued fraction is expressed by integrals, is even more remarkable, because its validity is not obvious at all. Although the case, in which $a=2$, was already found before and its value was given by the quadrature of the circle, the remaining cases cannot be derived from the same principles. Hence if this continued fraction in converted into a series in the way described at the beginning of this paper, one gets to such an intricate series that its sum cannot  be calculated by any means, except in the case $a=2$. Therefore, I have already worked quite hard before to demonstrate the validity of this theorem and  find a way to get to this continued fraction a priori; this investigation seemed the more difficult to me, the more useful I believed it to be for other investigations. But as long as all effort was without success concerning this task, I highly regretted that the method used by \textsc{Brouncker} was never explained and is probably completely lost.

\paragraph*{§18}

According to \textsc{Wallis's} \textsc{Brouncker} was led to this form by interpolation of this series

\[
\frac{1}{2}+\frac{1\cdot 3}{2 \cdot 4}+\cfrac{1 \cdot 3 \cdot 5}{2 \cdot 4 \cdot 6}+\text{etc.};
\]
for, \textsc{Wallis} had demonstrated that the intermediate terms of this series  give the quadrature of the circle. And even the beginning of this interpolation, also considered by \textsc{Brouncker} himself, was given in his book. Then it is said he\footnote{By this Euler means  Brouncker. In general Euler describes what he read in Wallis's book, what Wallis thought, how Brouncker found this continued fraction.} to resolve the single fractions $\frac{1}{2}$, $\frac{3}{4}$, $\frac{5}{6}$ etc. into two factors, which all constitute a continuous progression. So, if it was

\[
AB=\frac{1}{2}, \quad CD=\frac{2}{3}, \quad EF=\frac{5}{6}, \quad GH=\frac{7}{8} \quad \text{etc.}
\]
and the quantities $A$, $B$, $C$, $D$, $E$ etc. constitute a continuous progression that series becomes this one

\[
AB+ABCD+ABCDEF+\text{etc.};
\] 
and reduced to this form  it is easily interpolated.  Hence the term, whose index is $\frac{1}{2}$, is $=A$ and the term corresponding to the index $\frac{3}{2}$ is $=ABC$ and so forth. Hence this whole interpolation is reduced to the resolution of the single fractions into two factors each.

\paragraph*{§19}

But from the law of continuity it will be
\[
BC=\frac{2}{3}, \quad DE=\frac{4}{5}, \quad FG=\frac{6}{7} \quad \text{etc.}
\]
Therefore, because it is

\[
A=\frac{1}{2B}, \quad B=\frac{2}{3C}, \quad C=\frac{3}{4D}, \quad D=\frac{4}{5E} \quad \text{etc.},
\]
one immediately obtains

\[
A=\frac{1 \cdot 3 \cdot 3 \cdot 5 \cdot 5  \cdot 7}{2 \cdot 2 \cdot 4 \cdot 4 \cdot 6 \cdot 6} \quad \text{etc.},
\]
which is indeed the formula given by \textsc{Wallis} first, by which he expressed the quadrature of the circle, which formula differs from the  expression given by \textsc{Brouncker} a lot. Therefore, because this formula  is found so easily by investigating the interpolation this way, one has to wonder even more that \textsc{Brouncker} having used the same approach got to such a different expression; for, there seems to be no other way leading that continued fraction. And it is indeed not very likely that \textsc{Brouncker} wanted to express the value $A$ by a continued fraction on purpose, but rather following a certain peculiar method discovered it by accident, because at that time continued fractions were completely unknown and were introduced  on this occasion for the first time. From this it is possible to conclude that there is an obvious method leading to those continued fractions, even though it seems to be rather mysterious now.

\paragraph*{§20}

Although I spend  a lot of time on efforts to rediscover this method without any success, I nevertheless found another way to interpolate  series of this kind by continued fractions; but it led to expressions very different from the \textsc{Brounckerian} counterparts. Nevertheless, I think that it will be useful to explain this method, because applying it one finds continued fractions, whose values are known from elsewhere and can be exhibited by quadratures. For, then I will explain another method, to express the values of any continued fraction by quadratures. And using the result found by both methods  extraordinary comparisons of integral formulas will result, at least in that case, in which a definite value is attributed to the variable after the integration; I exhibited many comparisons of such kind in the preceding dissertation\footnote{Here, Euler again refers to his paper "De productis ex infinitis factoribus ortis".} on infinite products consisting of constant factors.

\paragraph*{§21}

To explain this method of interpolation I found, let this very general series be propounded 

\[
\frac{p}{p+2q}+\frac{p(p+2r)}{(p+2q)(p+2q+2r)}++\frac{p(p+2r)(p+4r)}{(p+2q)(p+2q+2r)(p+2q+4r)}+\text{etc.},
\]
whose term corresponding to the index $\frac{1}{2}$ is $=A$, the term corresponding to the index $\frac{3}{2}$ is $=ABC$, the term corresponding to the index $\frac{5}{2}$ is $=ABCDE$ etc. Hence it will be

\[
AB=\frac{p}{p+2q}, \quad CD=\frac{p+2r}{p+2q+2r}, \quad EF=\frac{p+4r}{p+2q+4r} \quad \text{etc.}
\]
and from the law of continuity

\[
BC=\frac{p+r}{p+2q+r}, \quad DE=\frac{p+3r}{p+2q+3r}, \quad FG=\frac{p+5r}{p+2q+5r}
\]
and so on.

\paragraph*{§22}

To get rid of the fractions put

\[
A=\frac{a}{p+2q-r}, \quad B=\frac{b}{p+2q}, \quad C=\frac{c}{p+2q+r}, \quad D=\frac{d}{p+2q+2r} \quad \text{etc.}
\]
and it will be

\[
ab=(p+2q-r)p, \quad bc=(p+2q)(p+r), \quad cd=(p+2q+r)(p+2r),
\]
\[
de=(p+2q+2r)(p+3r) \quad \text{etc.}
\]
Now let be

\[
a=m-r+\frac{1}{\alpha}, \quad  b=m+\frac{1}{\beta}, \quad c=m+r+\frac{1}{\gamma}, \quad
\]
\[
d=m-2r+\frac{1}{\delta}, \quad e=m+3r+\frac{1}{\epsilon} \quad \text{etc.},
\]
in which substitutions the integer parts constitute an arithmetic progression, whose constant difference is $r$; this is actually necessary because of the progression of those factors itself. Therefore, having substituted these values and, for the sake of brevity, having put

\[
p^2+2pq-pr-m^2+mr=P
\]
and

\[
2r(p+q-m)=Q
\]
the following equations will result

\begin{alignat*}{19}
&  && && &&P && \alpha \beta &&- && && &&(m&&-&&r &&)&&\alpha &&= && && &&m&&\beta &&+1,\\
& (P&&+&&&&Q)&& \beta \gamma  &&- && && && && && && m &&\beta &&= &&(m&&+&&r)&&\gamma &&+1,\\
& (P&&+&&2&&Q)&& \gamma \delta &&- && && &&(m&&+&& &&r)&&\gamma &&= &&(m&&+&&2r)&&\delta &&+1,\\
& (P&&+&&3&&Q)&&  \delta \varepsilon &&- && && &&(m&&+&&2&&r)&&\delta &&= &&(m&&+&&3r)&&\varepsilon &&+1,\\
& && && && && && && && && &&&&  &&&&&&\text{etc.}
\end{alignat*}

\paragraph*{§23}

Therefore, using these equations the following comparisons of the letters $\alpha$, $\beta$, $\gamma$, $\delta$ etc. result:

\begin{alignat*}{9}
&\alpha &&=\frac{m\beta +1}{P\beta-(m-r)}=\frac{m}{P}+\cfrac{p(p+2q-r):P^2}{-\cfrac{m-r}{P}+\beta}\\
&\beta &&=\frac{(m+r)\gamma +1}{(P+Q)\gamma -m}=\frac{m+r}{P+Q}+\cfrac{(p+r)(p+2q):(P+Q)^2}{-\cfrac{m}{P+Q}+\gamma}\\
&\gamma &&=\frac{(m+2r)\delta +1}{(P+2Q)\delta -(m+r)}=\frac{m+2r}{P+2Q}+\cfrac{(p+2r)(p+2q+r):(P+2Q)^2}{-\cfrac{m-r}{P+2Q}+\delta}\\
& && \qquad \qquad \qquad \qquad \qquad \qquad \qquad \text{etc.}
\end{alignat*}
So, if for the sake of brevity one puts

\[
p^2+2pq-mp-mq+qr=R
\]
and

\[
pr+qr-mr=S
\]
and the values of the assumed letters are continuously substituted into each other, the following continued fraction will result

\[
\alpha = \frac{m}{P}+\cfrac{p(p+2q-r):P^2}{\cfrac{2rR}{P(P+Q)}+\cfrac{(p+r)(p+2q):(P+Q)^2}{\cfrac{2r(R+S)}{(P+Q)(P+2Q)}+\cfrac{(p+2r)(p+2q+r):(P+2Q)^2}{\cfrac{2r(R+2S)}{(P+2Q)(P+3Q)}+\text{etc.}}}}
\]

\paragraph*{§24}

But because it is $a=m-r+\frac{1}{\alpha}$, one will have

\[
a=m-r+\cfrac{P}{m+\cfrac{p(p+2q-r)(P+Q)}{2rR+\cfrac{(p+r)(p+2q)P(P+2Q)}{2r(R+S)+\cfrac{(p+2r)(p+2q+r)(P+Q)(P+3Q)}{2r(R+2S)+\text{etc.}}}}}
\]
Hence the term of the propounded series

\[
\frac{p}{p+2q}+\frac{p(p+2r)}{(p+2q)(p+2q+2r)}+\frac{p(p+2r)(p+4r)}{(p+2q)(p+2q+2r)(p+2q+4r)}+\text{etc.}
\]
corresponding to the  index $\frac{1}{2}$ will be

\[
=A=\frac{a}{p+2q-r}.
\]
But because the general term  corresponding to the index $n$ of this series is

\[
=\frac{\int y^{p+2q-1}dy(1-y^{2r})^{n-1}}{\int y^{p-1}dy(1-y^{2r})^{n-1}}
\]
the  continued fraction found or the value of the letter $a$ is

\[
=(p+2q-r)\frac{\int y^{p+2q-1}dy:\sqrt{1-y^{2r}}}{\int y^{p-1}dy:\sqrt{1-y^{2r}}},
\]
 having put $y=1$ after both integrations, of course.
\paragraph*{§25}

But because the arbitrary letter $m$ is contained in our continued fraction, one will have innumerable continued fractions, whose value is the same and  hence known; therefore, it will be helpful, to consider the principal continued fractions. Hence at first let be

\[
m-r=p\quad \text{or} \quad m=p+r;
\]
it will be 

\[
P=2p(q-r), \quad Q=2r(q-r), \quad R=p(q-r) \quad \text{and} \quad S=r(q-r),
\]
whence it will be

\[
a=p+\cfrac{2p(q-r)}{p+r+\cfrac{(p+2q-r)(p+r)}{r+\cfrac{(p+2q)(p+2r)}{r+\cfrac{(p+2q+r)(p+3r)}{r+\text{etc.}}}}}
\]
But if $r>q$ so that the fraction does not become negative, it will be

\[
a=\cfrac{p}{1+\cfrac{2(r-q)}{p+2q-r+\cfrac{(p+2q-r)(p+r)}{r+\cfrac{(p+2q)(p+2r)}{r+\cfrac{(p+2q+r)(p+3r)}{r+\text{etc.}}}}}}
\]
\paragraph*{§26}

Now let $m$ be $=p+q$, in which case $Q$ and $S$ vanish; but it will be 

\[
P=q(r-q) \quad \text{and} \quad R=q(r-q)
\]
and hence it will result

\[
a=p+q-r+\cfrac{q(r-q)}{p+q+\cfrac{p(p+2q-r)}{2r+\cfrac{(p+r)(p+2q)}{2r+\cfrac{(p+2r)(p+2q+r)}{2r+\text{etc.}}}}},
\]
which continued fraction is therefore equal to the preceding ones, even though the forms  are different.

\paragraph*{§27}

Now put $m=p+2q$ and it will be

\begin{alignat*}{9}
&P&&=2q(r-p-2q)=-2q(p+2q-r),\\
&Q&&=-2qr,\\
&R&&=-q(p+2q-r)
\intertext{and}
&S&&=-qr.
\end{alignat*}
Therefore, one will obtain the following continued fraction:

\[
a=p+2q-r-\cfrac{2q(p+2q-r)}{p+2q+\cfrac{p(p+2q)}{r+\cfrac{(p+r)(p+2q+r)}{r+\cfrac{(p+2r)(p+2q+2r)}{r+\text{etc.}}}}}
\]
Hence innumerable continued fractions result and their  value is the same, namely $a$, which value  was found by means of integral formulas to be

\[
=(p+2q-r)\frac{\int y^{p+2q-1}dy:\sqrt{1-y^{2r}}}{\int y^{p-1}dy:\sqrt{1-y^{2r}}}=(p+2q-2r)\frac{\int y^{p+2q-2r-1}dy:\sqrt{1-y^{2r}}}{\int y^{p-1}dy:\sqrt{1-y^{2r}}}.
\]

\paragraph*{§28}

Before we proceed, let us consider some special cases. Let $r=2q$ and it will be

\[
a=p\frac{\int y^{p+2q-1}:\sqrt{1-y^{4q}}}{\int y^{p-1}dy:\sqrt{1-y^{4q}}}.
\]
Therefore, because it is

\begin{alignat*}{9}
&P&&=p^2+2mq-m^2,\\
&Q&&=4q(p+q-m),\\
&R&&=p^2+2pq+2qq-mp-mq
\intertext{and}
&S&&=2q(p+q-m),
\end{alignat*}
it will be in general

\[
a=m-2q+\cfrac{P}{m+\cfrac{p^2(P+Q)}{4qR+\cfrac{(p+2q)^2P(P+2Q)}{4q(R+S)+\cfrac{(p+4q)^2(P+Q)(P+3Q)}{4q(R+2S)+\text{etc.}}}}}
\]

\paragraph*{§29}

But if we substitute those different values for $m$, the following determined continued fraction will result

\[
a=p-\cfrac{2pq}{p+2q+\cfrac{p(p+2q)}{2q+\cfrac{(p+2q)(p+4q)}{2q+\cfrac{(p+4q)(p+6q)}{2q+\text{etc.}}}}}
\]
or instead of this continued fraction because of $r>q$

\[
a=\cfrac{p}{1+\cfrac{2q}{p+\cfrac{p(p+2q)}{2q+\cfrac{(p+2q)(p+4q)}{2q+\cfrac{(p+4q)(p+6q)}{2q+\text{etc.}}}}}}
\]
Further, from § $26$ one obtains this fraction for this case

\[
a=p-q+\cfrac{qq}{p+q+\cfrac{pp}{4q+\cfrac{(p+2q)^2}{4q+\cfrac{(p+4q)^2}{4q+\text{etc.}}}}}
\]
But finally, § $27$ will give this continued fraction

\[
a=p-\frac{2pq}{p+2q+\cfrac{p(p+2q)}{2q+\cfrac{(p+2q)(p+4q)}{2q+\cfrac{(p+4q)(p+6q)}{2q+\text{etc.}}}}}
\]
which agrees to the one exhibited first here, so that one has only two simpler continued fractions for this case $r=2q$.

\paragraph*{§30}

Now further put $q=p=1$ that it is

\[
a= \frac{\int yy dy :\sqrt{1-y^4}}{\int dy :\sqrt{1-y^4}};
\]
first, it will be

\[
a=-\cfrac{2}{3+\cfrac{1 \cdot 3}{2+\cfrac{3 \cdot 5}{2+\cfrac{5 \cdot 7}{2+\text{etc.}}}}}
\]
Then one will have

\[
a=\cfrac{1}{2+\cfrac{1}{4+\cfrac{9}{4+\cfrac{25}{4+\cfrac{49}{4+\text{etc.}}}}}}
\]
Hence it follows that it will be

\[
\frac{\int dy: \sqrt{1-y^4}}{\int yydy:\sqrt{1-y^4}}=2+\cfrac{1}{4+\cfrac{9}{4+\cfrac{25}{4+\cfrac{49}{4+\text{etc.}}}}}
\]
which case is contained in the expression given in § $16$; hence that formula, not rigorously proven yet, is, at least, confirmed a bit more. Then having put $a=3$ there, it will be

\[
3 \frac{\int x^4 dx :\sqrt{1-x^4}}{\int xxdx:\sqrt{1-x^4}}=\frac{\int dx :\sqrt{1-x^4}}{\int xxdx:\sqrt{1-x^4}}=2+\cfrac{1}{4+\cfrac{9}{4+\cfrac{25}{4+\cfrac{49}{4+\text{etc.}}}}}
\]
so that it is now certain that the formula exhibited in § $16$ is true in the cases, in which it is either $a=2$ or $a=3$; but it will soon be shown to be true in all cases.

\paragraph*{§31}

Let $q=\frac{1}{2}$ and $p=1$; while it still is $r=2q=1$, it will be

\[
a= \frac{\int ydy:\sqrt{1-y^2}}{\int dy:\sqrt{1-y^2}}=\frac{2}{\pi},
\]
where $\pi$ denotes the circumference of the circle, whose diameter is $=1$. Therefore, it will be in general

\[
P=1+m-m^2, \quad Q=2-2m,
\]
\[
R=\frac{5-3m}{2} \quad \text{and} \quad S=\frac{3-2m}{2}
\]
and therefore

\[
a=m-1+\cfrac{1+m-m^2}{m+\cfrac{1^2(4-m+m^2)}{5-3m+\cfrac{2^2(1+m-m^2)(7-3m-m^2)}{8-5m+\cfrac{3^2(4-m-m^2)(10-5m-m^2)}{11-7m+\text{etc.}}}}}
\]
But in the special cases expanded above it will be

\begin{alignat*}{9}
&\frac{\pi}{2}&&=\cfrac{1}{1-\cfrac{1}{2+\cfrac{1 \cdot 2}{1+\cfrac{2 \cdot 3}{1+\cfrac{3 \cdot 4}{1+\text{etc.}}}}}} \quad &&=1&&+\cfrac{1}{1+\cfrac{1 \cdot 2}{1+\cfrac{2 \cdot 3}{1+\cfrac{3 \cdot 4}{1+\text{etc.}}}}}
\intertext{and}
&\frac{\pi}{2}&&=\cfrac{1}{{\frac{1}{2}+\cfrac{1 : 4}{\frac{3}{2}+\cfrac{1^2}{2+\cfrac{2^2}{2+\cfrac{3^2}{2+\text{etc.}}}}}}} \quad &&=2&&-\cfrac{1}{2+\cfrac{1^2}{2+\cfrac{2^2}{2+\cfrac{3^2}{2+\text{etc.}}}}}
\end{alignat*}

\paragraph*{§32}

To understand the use of these formulas for interpolations, let this series be propounded
\[
\frac{2}{1}+\cfrac{2 \cdot 4}{1 \cdot 3}+\cfrac{2 \cdot 4 \cdot 6}{1 \cdot 3 \cdot 5}+\text{etc.}
\]
whose term corresponding to the index $\frac{1}{2}$ has to be found; set this term $=A$; therefore, it will be

\[
p=2, \quad r=1 \quad \text{and} \quad q=-\frac{1}{2}.
\]
Now put

\[
A=\frac{a}{p+2q-r};
\]
it will be

\[
A=\frac{a}{0},
\]
whence the inconvenience of the given formulas, if $p+2q-r=0$, is seen clearly.  It nevertheless is possible to complete this task by finding the term corresponding to the index $\frac{3}{2}$; if this term was $=Z$, it will be $A=\frac{2}{3}Z$; but $\frac{1}{2}Z$ will be the term corresponding to the index $\frac{1}{2}$ of this series

\[
\frac{4}{3}+\frac{4 \cdot 6}{3 \cdot 5}+\cfrac{4 \cdot 6 \cdot 8}{3 \cdot 5 \cdot 7}+\text{etc.},
\] 
which compared to the general series gives

\[
p=4, \quad r=1, \quad q=-\frac{1}{2},
\]
so that it is

\[
Z=\frac{2 \int y^2dy:\sqrt{1-y^2}}{\int y^3 dy:\sqrt{1-y^2}}=\frac{3 \int dy:\sqrt{1-y^2}}{2 \int dy:\sqrt{1-y^2}}=\frac{3}{4}\pi
\]
and $A=\frac{\pi}{2}$. But because by § $24$ it is

\[
Z=a \quad \text{and} \quad A=\frac{2}{3}Z=\frac{2}{3}a,
\]
at first, because of

\[
P=8+m-m^2, \quad Q=7-2m
\]
\[
R=\frac{23-7m}{2} \quad \text{and} \quad S=\frac{7-2m}{2} 
\]
it will be in general

\[
A=\frac{3}{2}a=\frac{\pi}{2}
\]
\[
=\frac{2(m-1)}{3}+\cfrac{2(8+m-m^2)}{3m+\cfrac{2 \cdot 4 \cdot 3(15-m-m^2)}{23-7m+\cfrac{3 \cdot 5(8+m-m^2)(22-3m-m^2)}{30-9m+\cfrac{4 \cdot 6(15-m-m^2)(29-5m-^2)}{37-11m+\text{etc.}}}}}
\]

\paragraph*{§33}

But by expanding the particular cases it will be

\[
a=\frac{3}{4}\pi=4-\cfrac{12}{5+\cfrac{2 \cdot 5}{1+\cfrac{3 \cdot 6}{1+\cfrac{4 \cdot 7}{1+\text{etc.}}}}}=\cfrac{4}{1+\cfrac{1 \cdot 3}{2+\cfrac{2 \cdot 5}{1+\cfrac{3 \cdot 6}{1+\cfrac{4 \cdot 7}{1+\text{etc.}}}}}}
\]
or also

\[
\frac{3}{4}\pi=1+\cfrac{3}{1+\cfrac{1 \cdot 4}{1+\cfrac{2 \cdot 5}{1+\cfrac{3 \cdot 6}{1+\cfrac{4 \cdot 7}{1+\text{etc.}}}}}}
\]
In like manner by § $26$ one will find

\[
a=\frac{3}{4}\pi= \frac{5}{2}-\cfrac{3:4}{\frac{7}{2}+\cfrac{2 \cdot 4}{2+\cfrac{3 \cdot 5}{2+\cfrac{4 \cdot 6}{2+\cfrac{5 \cdot 7}{2+\text{etc.}}}}}}=2+\cfrac{1}{2+\cfrac{1 \cdot 3}{2+\cfrac{2 \cdot 4}{2+ \cfrac{3 \cdot 5}{2+\cfrac{4 \cdot 6}{2+\text{etc.}}}}}}
\]
Finally, the case considered in § $27$ will give

\[
a= \frac{3}{4}\pi = 2+\cfrac{2}{3+\cfrac{3 \cdot 4}{1+\cfrac{4 \cdot 5}{1+\cfrac{5 \cdot 6}{1+\text{etc.}}}}}
\]
or

\[
\frac{\pi}{2}=1+\cfrac{1}{1+\cfrac{1 \cdot 2}{1+\cfrac{2 \cdot 3}{1+\cfrac{3 \cdot 4}{1+\cfrac{4 \cdot 5}{1+\text{etc.}}}}}}
\]
which expression agrees with the one exhibited in § $31$.

\paragraph*{§34}

Hence we obtained innumerable continued fractions from this method of interpolation; and  the values of these fractions can be assigned by quadratures of curves or integral formulas. But because these continued fractions are irregular at the beginning, just truncate them in such a way, that the continued fractions proceed regularly. Therefore, using § $25$ and putting

\[
p+2q-r=f \quad \text{and} \quad p+r=h 
\]
the following equation will result

\[
r+\cfrac{fh}{r+\cfrac{(f+r)(h+r)}{r+\cfrac{(f+2r)(h+2r)}{r+\text{etc.}}}}
\]
\[
=\frac{h(f-r)\int y^{h+r-1}dy:\sqrt{1-y^{2r}}-f(h-r)\int y^{f+r-1}dy: \sqrt{1-y^{2r}}}{f \int y^{f+t-1}dy :\sqrt{1-y^{2r}}-h\int y^{h+r-1} dy :\sqrt{1-y^{2r}}},
\]
which equation is always correct,  except for the case $f=h$. But in order to treat the case, in which $f=h$, just put $f=h+dw$ and one will find

\[
\frac{\int y^{h+r+dw-1}dy:\sqrt{1-y^{2r}}}{\int y^{h+r-1}dy :\sqrt{1-y^{2r}}}=1-rdw \int \frac{dx}{x^{r+1}}\int \frac{x^{h+2r-1}dx}{1-x^{2r}}
\]
after having put $x=1$ after the integration. Hence it will be

\[
r+\cfrac{hh}{r+\cfrac{(h+r)^2}{r+\cfrac{(h+2r)^2}{r+\text{etc.}}}}
\]
\[
= \frac{r+hr(h-r) \int \frac{dx}{x^{r+1}}\int \frac{x^{h+2r-1}dx}{1-x^{2r}}}{1-hr\int \frac{dx}{x^{r+1}}\int \frac{x^{h+2r-1}dx}{1-x^{2r}}}= \frac{r(h-r)^2\int \frac{dx}{x^{r+1}}\int \frac{x^{h-1}dx}{1-x^{2r}}}{1-r(h-r) \int \frac{dx}{x^{r+1}}\int \frac{x^{h-1}dx}{1-x^{2r}}}.
\]
But from the nature of the integrals one has

\[
\int \frac{dx}{x^{r+1}}\int \frac{x^{h+2r-1}dx}{1-x^{2r}}=\frac{-1}{rx^r}\int \frac{x^{h+2r-1}dx}{1-x^{2r}}+\frac{1}{r}\int \frac{x^{h+r-1}dx}{1-x^{2r}}=\frac{1}{r}\int \frac{x^{h+r-1}dx}{1+x^r}
\]
having put $x=1$. Therefore, one will have

\[
r+\cfrac{hh}{r+\cfrac{(h+r)^2}{r+\text{etc.}}}=\frac{r+h(h-r)\int \frac{x^{h+r-1}dx}{1+x^r}}{1-h \int \frac{x^{h+r-1}dx}{1+x^r}}=\frac{1-(h-r) \int \frac{x^{h-1}dx}{1+x^r}}{\int \frac{x^{h-1}dx}{1+x^r}};
\]
this formula agrees with the one  given in § $7$.

\paragraph*{§35}

In like manner  it follows from § $26$ by putting $p=f$ and $p+2q-r=h$ that it will be

\[
2r+\cfrac{fh}{2r+\cfrac{(f+r)(h+r)}{2r+\cfrac{(f+2r)(h+2r)}{2r+\text{etc.}}}}
\]
\[
=\frac{2(r-f)(r-h) \int \frac{y^{f-1}dy}{\sqrt{1-y^{2r}}}-h(f+h-3r)\int \frac{y^{h+r-1}dy}{\sqrt{1-y^{2r}}}}{2h \int \frac{y^{h+r-1}dy}{\sqrt{1-y^{2r}}}-(f+h-r)\int \frac{y^{f-1}dy}{\sqrt{1-y^{2r}}}}
\]
But because this formula  remains the same, if $f$ and $h$ are commuted, it is manifest that is has to be

\[
\frac{h \int y^{h+r-1}dx:\sqrt{1-y^{2r}}}{y^{f-1}dy: \sqrt{1-y^{2r}}}= \frac{f \int y^{f+r-1}dy:\sqrt{1-y^{2r}}}{\int y^{h-1}dy:\sqrt{1-y^{2r}}}
\]
after having put $y=1$ after every integration. But this theorem is already contained  in those I exhibited in the preceding dissertation\footnote{Euler refers to his paper E122 again.} on products consisting of infinitely many factors; for, there I produced many theorems of this kind and proved them.

\paragraph*{§36}

But here in like manner the case, in which it $f=h+r$, is worth one's attention; for, in this case so the numerator as the denominator of the found continued fraction vanishes. But having put $f=h+r+dw$ as before and having done the calculation it will result

\[
2r+\cfrac{h(h+r)}{2r+\cfrac{(h+r)(h+2r)}{2r+\cfrac{(h+2r)(h+3r)}{2r+\text{etc.}}}}= \frac{h+2h(r-h)\int \frac{x^{h-r}dx}{1+x^r}}{-1+2h\int \frac{x^{h-r}dx}{1+x^r}}.
\]
Hence if one puts $h=r=1$, one will have

\[
2+\cfrac{1 \cdot 2}{2+\cfrac{2 \cdot 3}{2+\cfrac{3 \cdot 4}{2+\cfrac{4 \cdot 5}{2+\text{etc.}}}}}=\frac{1}{2 \log 2-1}.
\]
If additionally the equation in § $27$ is treated in the same way, a form very similar to that one I found § $25$, will result.

\paragraph*{§37}

Having explained how the interpolation of series is reduced to continued fractions, I return to \textsc{Brounckerian} expressions and will give a natural method not only to get to those, but also a method probably used by \textsc{Brouncker} himself. But the continued fractions found up to this point differ from the  ones given by \textsc{Brouncker} a lot, because the values of the letters $A$, $B$, $C$, $D$ etc., using  the method explained, depend on each other in such a way that they can  easily be compared to each other, but using the \textsc{Brouckerian} method they turned out to be so different to each other such that their relation is not seen. This difference itself finally led me to the invention of another method I will explain now.

\paragraph*{§38}

But before I explain this method of interpolation, it will be convenient to mention the following very general lemma in advance. If there were innumerable quantities $\alpha$, $\beta$, $\gamma$, $\delta$, $\varepsilon$ etc., which depend on each other in such a way, that it is

\begin{alignat*}{9}
&\alpha \beta &&-m \alpha &&-n \beta &&-\varkappa  &&=0,\\
&\beta \gamma &&-(m+s)\beta &&-(n+s)\gamma &&-\varkappa &&=0,\\
& \gamma \delta &&-(m+2s)\gamma &&-(n+2s)\delta &&-\varkappa &&=0,\\
& \delta \varepsilon &&-(m+3s)\delta &&-(n+3s)\varepsilon &&-\varkappa &&=0\\
& &&  && \text{etc.},
\end{alignat*}
and  the letters $\alpha$, $\beta$, $\gamma$, $\delta$ etc. obtain the following values

\begin{alignat*}{9}
&\alpha &&=m &&+n &&-s&&-\frac{ss-ms+ns+\varkappa}{a},\\
&\beta &&=m &&+n &&+s&&-\frac{ss-ms+ns+\varkappa}{b},\\
&\gamma &&=m &&+n &&+3s&&-\frac{ss-ms+ns+\varkappa}{c},\\
&\delta &&=m &&+n &&+5s&&-\frac{ss-ms+ns+\varkappa}{d}\\
& && && && && \text{etc.}
\end{alignat*}
the upper equations will be transformed into the following similar ones

\begin{alignat*}{9}
&ab&&-(m-s)a &&-(n+s)b&&-ss&&+ms&&-ns&&-\varkappa &&=0,\\
&bc&&-mb &&-(n+2s)c&&-ss&&+ms&&-ns&&-\varkappa &&=0,\\
&cd&&-(m+s)c &&-(n+3s)d&&-ss&&+ms&&-ns&&-\varkappa &&=0,\\
&de&&-(m+2s)d &&-(n+4s)4&&-ss&&+ms&&-ns&&-\varkappa &&=0\\
& && && &&  \text{etc.}
\end{alignat*}
This fact itself that formulas of the same structure arose was the origin of those substitutions.

\paragraph*{§39}

If now in like manner these last equations are transformed  into similar ones by means of suitable substitutions, one will find the following substitutions for $a$, $b$, $c$, $d$ etc.

\begin{alignat*}{9}
&a &&=m+n&&-s &&+\frac{4ss-2ms+2ns+\varkappa}{a_1},\\
&b &&=m+n&&+s &&+\frac{4ss-2ms+2ns+\varkappa}{b_1},\\
&c &&=m+n&&+3s &&+\frac{4ss-2ms+2ns+\varkappa}{c_1},\\
&d &&=m+n&&+5s &&+\frac{4ss-2ms+2ns+\varkappa}{d_1}\\
& && && &&\text{etc.}
\end{alignat*}
having substituted these values the following equations will result

\begin{alignat*}{9}
&a_1b_1&&-(m-2s)&&a_1&&-(n+2s)b_1&&-4ss+2ms-2ns-\varkappa &&=0,\\
&b_1c_1&&-(m-s)&&b_1&&-(n+3s)c_1&&-4ss+2ms-2ns-\varkappa &&=0,\\
&c_1d_1&&-m&&c_1&&-(n+4s)d_1&&-4ss+2ms-2ns-\varkappa &&=0,\\
&d_1e_1&&-(m+s)&&d_1&&-(n+5s)e_1&&-4ss+2ms-2ns-\varkappa &&=0\\
& && && && && \text{etc.}
\end{alignat*}

\paragraph*{§40}

Proceeding further one will have to put

\begin{alignat*}{9}
&a_1&&=m+n&&-s&&+\frac{9ss-3ms-3ns+\varkappa}{a_2},\\
&b_1&&=m+n&&+s&&+\frac{9ss-3ms-3ns+\varkappa}{b_2},\\
&c_1&&=m+n&&+3s&&+\frac{9ss-3ms-3ns+\varkappa}{c_2},\\
& && && && \text{etc.}
\end{alignat*}
And from those substitutions these equations result

\begin{alignat*}{9}
&a_2b_2&&-(m-3s)&&-(n+3s)b_2&&-9ss+3ms-3ns-\varkappa &&=0,\\
&b_2c_2&&-(m-2s)&&-(n+4s)c_2&&-9ss+3ms-3ns-\varkappa &&=0,\\
&c_2d_2&&-(m-s)&&-(n+5s)d_2&&-9ss+3ms-3ns-\varkappa &&=0\\
& && && &&\text{etc.}
\end{alignat*}

\paragraph*{§41}

If these substitutions are now continued to infinity and the following values are always substituted in the preceding ones, the values of the letters $\alpha$, $\beta$, $\gamma$, $\delta$ etc. are expressed by the following continued fractions

\begin{alignat*}{9}
&\alpha &&=m+n&&-s&&+\cfrac{ss-ms+ns-\varkappa}{m+n-s+\cfrac{4ss-2ms+2ns+\varkappa}{m+n-s+\cfrac{9ss-3ms+3ns+\varkappa}{m+n-s+\cfrac{16ss-4ms+4ns+\varkappa}{m+n-s+\text{etc.}}}}}\\
&\beta &&=m+n&&+s&&+\cfrac{ss-ms+ns-\varkappa}{m+n+s+\cfrac{4ss-2ms+2ns+\varkappa}{m+n+s+\cfrac{9ss-3ms+3ns+\varkappa}{m+n+s+\cfrac{16ss-4ms+4ns+\varkappa}{m+n+s+\text{etc.}}}}}\\
&\gamma &&=m+n&&+3s&&+\cfrac{ss-ms+ns-\varkappa}{m+n+3s+\cfrac{4ss-2ms+2ns+\varkappa}{m+n+3s+\cfrac{9ss-3ms+3ns+\varkappa}{m+n+3s+\cfrac{16ss-4ms+4ns+\varkappa}{m+n+3s+\text{etc.}}}}}\\
& && && &&\text{etc.}
\end{alignat*} 
which continued fractions are similar enough to those \textsc{Brouncker} gave, while the following are not contained in the preceding ones.

\paragraph*{§42}

But that the use of these formulas for interpolations becomes more apparent, let this series be propounded

\[
\frac{p}{p+2q}+\frac{p(p+2r)}{(p+2q)(p+2q+2r)}+\frac{p(p+2r)((p+4r)}{(p+2q)(p+2q+2r)(p+2q+4r)}+\text{etc.}
\] 
and let the term corresponding to the index $\frac{1}{2}$  be $=A$, the term  corresponding to the index $\frac{3}{2}=ABC$, the term corresponding to the index $\frac{5}{2}=ABCDE$ and so on. Hence it will be

\[
AB=\frac{p}{p+2q}, \quad CD=\frac{p+2r}{p+2q+2r}, \quad EF=\frac{p+4r}{p+2q+4r} \quad \text{etc.}
\]
and it will be

\begin{alignat*}{9}
&ab&&=p(p+2q-r),\\
&bc&&=(p+r)(p+2q),\\
&cd&&=(p+2r)(p+2q+r),\\
&de&&=(p+3r)(p+2q+2r)\\
& && \text{etc.}
\end{alignat*}
Now further let it be

\begin{alignat*}{9}
&a&&=p+q&&-r&&+\frac{g}{\alpha},\\
&b&&=p+q&& &&+\frac{g}{\beta},\\
&c&&=p+q&&+r&&+\frac{g}{\gamma},\\
&d&&=p+q&&+2r&&+\frac{g}{\delta}\\
&  && \text{etc.};
\end{alignat*}
having substituted and having put $g=q(r-q)$ these values the following equations will emerge

\begin{alignat*}{9}
&\alpha \beta &&(p+q-r)&&\alpha &&- (p+q)&&\beta &&-q(r-q)&&=0, \\
& \beta \gamma &&(p+q)&&\beta &&- (p+q+r)&&\gamma &&-q(r-q)&&=0, \\
&\gamma \delta &&(p+q+r)&& \gamma &&- (p+q+2r)&&\delta &&-q(r-q)&&=0, \\
&\delta \varepsilon &&(p+q+2r)&&\delta &&- (p+q+3r)&&\varepsilon &&-q(r-q)&&=0, \\
& \text{etc.}
\end{alignat*}
Having substituted all these values one will obtain the following continued fractions expressing the letters $a$, $b$, $c$, $d$ etc.:

\begin{alignat*}{9}
&a&&=p+q-r+\cfrac{qr-qq}{2(p+q-r)+\cfrac{2rr+qr-qq}{2(p+q-r)+\cfrac{6rr+qr-qq}{2(p+q-r)+\cfrac{12rr+qr-qq}{2(p+q-r)+\text{etc.}}}}}\\
&b&&=p+q+\cfrac{qr-qq}{2(p+q)+\cfrac{2rr+qr-qq}{2(p+q)+\cfrac{6rr+qr-qq}{2(p+q)+\cfrac{12rr+qr-qq}{2(p+q)+\text{etc.}}}}}\\
&c&&=p+q+r+\cfrac{qr-qq}{2(p+q+r)+\cfrac{2rr+qr-qq}{2(p+q+r)+\cfrac{6rr+qr-qq}{2(p+q+r)+\cfrac{12rr+qr-qq}{2(p+q+r)+\text{etc.}}}}}\\
&  &&\text{etc.}
\end{alignat*}

\paragraph*{§44}

But because the term of the propounded series corresponding to the index $n$, is

\[
=\frac{\int y^{p+2q-1}dy(1-y^{2r})^{n-1}}{\int y^{p-1}dy(1-y^{2r})^{n-1}},
\]
it will be

\[
A=\frac{q}{P+2q-r}= \frac{\int y^{p+2q-1}dy:\sqrt{1-y^{2r}}}{\int y^{p-1}dy :\sqrt{1-y^{2r}}}
\]
or
\[
a=(p+2q-r)\frac{\int y^{p+2q-1}dy:\sqrt{1-y^{2r}}}{\int y^{p-1}dy :\sqrt{1-y^{2r}}}.
\]
Then, because of $ab=p(p+2q-r)$, it will be

\[
b=\frac{b\int y^{p-1}dy:\sqrt{1-y^{2r}}}{\int y^{p+2q-1}dy :\sqrt{1-y^{2r}}}.
\]
But because by a theorem shown in the preceding dissertation\footnote{This means E122 again.} it is

\[
\frac{p\int y^{p-1}dy:\sqrt{1-y^{2r}}}{\int y^{f+r-1}dy :\sqrt{1-y^{2r}}}=\frac{f\int y^{f-1}dy:\sqrt{1-y^{2r}}}{\int y^{p+r-1}dy :\sqrt{1-y^{2r}}}=\frac{(f+r)\int y^{f+2r-1}dy:\sqrt{1-y^{2r}}}{\int y^{p+r-1}dy :\sqrt{1-y^{2r}}},
\]
just put $f=p+2q-r$; having done so it will be

\[
b=\frac{(p+2q)\int y^{p+2q+r-1}dy:\sqrt{1-y^{2r}}}{\int y^{p+r-1}dy :\sqrt{1-y^{2r}}}.
\]
Proceeding in the same way it will indeed be

\[
c=\frac{(p+2q+r)\int y^{p+2q+2r-1}dy:\sqrt{1-y^{2r}}}{\int y^{p+2r-1}dy :\sqrt{1-y^{2r}}}
\]
and

\[
d=\frac{(p+2q+2r)\int y^{p+2q+3r-1}dy:\sqrt{1-y^{2r}}}{\int y^{p+3r-1}dy :\sqrt{1-y^{2r}}}
\]
\[
\text{etc.}
\]

\paragraph*{§45}

Therefore, since it is known how these integral formulas proceed, one  can deduce that the value of this continued fraction

\[
p+q+mr+\cfrac{qr-qq}{2(p+q+mr)+\cfrac{2rr+qr-qq}{2(p+q+mr)+\cfrac{6rr+qr-qq}{2(p+q+mr)+\text{etc.}}}}
\]
is
\[
=(p+2q+mr)\frac{\int y^{p+2q+(m+1)r-1}dy:\sqrt{1-y^{2r}}}{\int y^{p+(m+1)r-1}dy :\sqrt{1-y^{2r}}}.
\]
Therefore, if one puts $p+q+mr=s$, so that $p=s-q-mr$, the following continued fraction will result

\[
s+\cfrac{qr-qq}{2s+\cfrac{2rr+qr-qq}{2s+\cfrac{6rr+qr-rr}{2s+\cfrac{12rr+qr-qq}{2s+\cfrac{20rr+qr-qq}{2s+\text{etc.}}}}}}
\]
whose value will therefore be this expression

\[
(q+s) \frac{\int y^{q+r+s-1}dy:\sqrt{1-y^{2r}}}{\int y^{r+s-q-1}dy:\sqrt{1-y^{2r}}}.
\]

\paragraph*{§46}

Because in like manner the value of this continued fraction

\[
s+r+\frac{qr-qq}{2(s+r)+\cfrac{2rr+qr-qq}{2(s+r)+\cfrac{6rr+qr-qq}{2(s+r)+\text{etc.}}}}
\]
is

\[
=(q+r+s)\frac{\int y^{p+2r+q-1}dy:\sqrt{1-y^{2r}}}{\int y^{p+2r-q-1}dy :\sqrt{1-y^{2r}}},
\]
the product of this two continued fractions will hence be

\[
=(s+q)(s+r-q),
\]
as the product of the integral formulas reveals. Hence  by the theorem given in the preceding dissertation\footnote{This also refers to E122.} it is

\[
\frac{f}{a}=\frac{\int x^{a-1}dx:\sqrt{1-x^{2r}} \cdot \int x^{a+r-1}dx:\sqrt{1-x^{2r}}}{\int x^{f-1}dx:\sqrt{1-x^{2r}} \cdot \int x^{f+r-1}dx:\sqrt{1-x^{2r}}};
\]
the product of integral formulas is easily reduced to this form.

\paragraph*{§47}

The  continued fraction found can be transformed into another more convenient form  such that the single numerators can be resolved into factors; so one will have this continued fraction

\[
s+\cfrac{q(r-q)}{2s+\cfrac{(r+q)(2r-q)}{2s+\cfrac{(2r+q)(3r-q)}{2s+\cfrac{(3r+q)(4r-q)}{2s+\text{etc.}}}}}
\]
whose value will therefore be

\[
=(q+s)\frac{\int y^{r+s+q-1}dy:\sqrt{1-y^{2r}}}{\int y^{r+s-q-1}dy:\sqrt{1-y^{2r}}}.
\]
Therefore, if $s$ is added to the continued fraction, so that the structure is the same everywhere, it will be

\[
\frac{(q+s)\int y^{r+s+q-1}dy:\sqrt{1-y^{2r}}+\int y^{r+s-q-1}dy:\sqrt{1-y^{2r}}}{\int y^{r+s+q-1}dy:\sqrt{1-y^{2r}}}
\]
\[
=2s+\cfrac{q(r-q)}{2s+\cfrac{(r+q)(2r-q)}{2s+\cfrac{(2r+q)(3r-q)}{2s+\cfrac{(3r+q)((4r-q)}{2s+\text{etc.}}}}}
\]

\paragraph*{§48}

If now one puts $r=2$ and $q=1$, all continued fractions exhibited by \textsc{Brouncker} will result; for, those fractions will all be contained in this continued fraction

\[
s+\cfrac{1}{2s+\cfrac{9}{2s+\cfrac{25}{2s+\cfrac{49}{2s+\cfrac{81}{2s+\text{etc.}}}}}}
\]
whose value will therefore be

\[
=(s+1)\frac{\int y^{s+2}dy:\sqrt{1-y^{2}}}{\int y^{s}dy:\sqrt{1-y^{2}}},
\]
which expression agrees perfectly with that one we gave above, before its truth was absolutely certain; see § $16$.

\paragraph*{§49}

Whereas up to now I gave many continued fractions, whose values can be expressed by integral formulas, I will now explain a direct method  how to get  vice versa from integral formulas to continued fractions. But this method is based on a reduction of one integral formula to two others, which reduction is quite similar to that usual one, by which the integration of a certain differential formula is reduced to the integration of another\footnote{By this Euler means integration by parts.}. Let the infinitely many integral formulas be of this kind

\[
\int Pdx, \quad \int PRdx, \quad \int PR^2dx, \quad \int PR^3dx, \quad \int PR^4dx \quad \text{etc.},
\]
and further let them be of such a nature that, if they  are integrated, they vanish for $x=0$, and then $x$ is put $=1$; hence it is as follows

\begin{alignat*}{9}
&~~~~a&&\int Pdx &&=~~~~b &&\int PRdx &&+~~~~c &&\int PR^2dx, \\
&(a+\alpha)&&\int PRdx &&=(b+\beta)&& \int PR^2dx &&+(c+\gamma) &&\int PR^3dx, \\
&(a+2\alpha)&&\int PR^2dx &&=(b+2\beta)&& \int PR^3dx &&+(c+2\gamma) &&\int PR^4dx, \\
&(a+3\alpha)&&\int PR^3dx &&=(b+3\beta)&& \int PR^4dx &&+(c+3\gamma) &&\int PR^5dx
\intertext{and in general}
&(a+n\alpha)&&\int PR^ndx &&=(b+n\beta)&& \int PR^{n+1}dx &&+(c+n\gamma) &&\int PR^{n+2}dx.
\end{alignat*}

\paragraph*{§50}

Therefore, if one has integral formulas of this kind one can easily construct a continued fraction from them. Because it is

\begin{alignat*}{9}
&\frac{\int Pdx}{\int PRdx}&&= ~~~\frac{b}{a}&&+\frac{c\int PR^2dx}{a \int PRdx},\\
&\frac{\int PRdx}{\int PR^2dx}&&= \frac{b+\beta}{a+\alpha}&&+\frac{(c+\gamma)\int PR^3dx}{(a+\alpha) \int PR^2dx},\\
&\frac{\int PR^2dx}{\int PR^3dx}&&= \frac{b+2\beta}{a+2\alpha}&&+\frac{(c+2\gamma)\int PR^4dx}{(a+2\alpha) \int PR^3dx},\\
&\frac{\int PR^3dx}{\int PR^4dx}&&= \frac{b+3\beta}{a+3\alpha}&&+\frac{(c+3\gamma)\int PR^5dx}{(a+3\alpha) \int PR^4dx}\\
& && && \text{etc.},
\end{alignat*}
by substituting each value in the preceding equation it will be

\[
\frac{\int Pdx}{\int PRdx}=\frac{b}{a}+\cfrac{c:a}{\cfrac{b+\beta}{a+\alpha}+\cfrac{(c+\gamma):(\alpha +a)}{\cfrac{b+2 \beta}{a+2 \alpha}+\cfrac{(c+2 \gamma):(a+2 \alpha)}{\cfrac{b+3 \beta}{a+3 \alpha}+\cfrac{(c+3 \gamma)(:(a+3 \alpha)}{\cfrac{b+4 \beta}{a+4 \alpha}+\text{etc.}}}}}
\]
But having inverted this expression and having cleared the partial fractions one will have this expression

\[
\frac{\int PRdx}{\int Pdx}=\cfrac{a}{b+\cfrac{(a+\alpha)c}{b+\beta+\cfrac{(a+2 \alpha)(c+\gamma)}{b+2 \beta+\cfrac{(a+3\alpha)(c+2 \gamma)}{b+3 \beta+\cfrac{(a+4 \alpha)(c+3 \gamma)}{b+4 \beta}+\text{etc.}}}}}
\]

\paragraph*{§51}

If, while $n$ was a negative integer, it was

\[
(a+n\alpha)\int PR^ndx=(b+n \beta) \int PR^{n+1}dx+(c+n \gamma) \int PR^{n+2}dx,
\]
one will have the following equations:

\begin{alignat*}{9}
&(a- \alpha)&&\int \frac{Pdx}{R}&&=(b-\beta)&&\int Pdx&&+(c- \gamma)\int PRdx,\\
&(a- 2\alpha)&&\int \frac{Pdx}{R^2}&&=(b-2\beta)&&\int \frac{Pdx}{R}&&+(c- 2\gamma)\int Pdx,\\
&(a- 3\alpha)&&\int \frac{Pdx}{R^3}&&=(b-3\beta)&&\int \frac{Pdx}{R^2}&&+(c- 3\gamma)\int \frac{Pdx}{R},\\
&(a- 4\alpha)&&\int \frac{Pdx}{R^4}&&=(b-4\beta)&&\int \frac{Pdx}{R^3}&&+(c- 3\gamma)\int \frac{Pdx}{R^2}\\
& && && && \text{etc.}
\end{alignat*}
Hence in like manner it is deduced

\begin{alignat*}{9}
&\frac{\int PRdx}{\int Pdx}&&=\frac{-(b-\beta)}{c-\gamma}&&+\frac{(a-\alpha)\int Pdx:R}{(c-\gamma)\int Pdx},\\
&\frac{\int Pdx}{\int Pdx:R}&&=\frac{-(b-2\beta)}{c-2\gamma}&&+\frac{(a-2\alpha)\int Pdx:R^2}{(c-2\gamma)\int Pdx:R},\\
&\frac{\int Pdx:R}{\int Pdx:R^2}&&=\frac{-(b-3\beta)}{c-3\gamma}&&+\frac{(a-3\alpha)\int Pdx:R^3}{(c-3\gamma)\int Pdx:R^2}\\
& && &&  \text{etc.}
\end{alignat*}
But from these equations one gets to

\[
\frac{\int Prdx}{Pdx}= \frac{-(b-\beta)}{c-\gamma}+\cfrac{(a-\alpha):(c-\gamma)}{\cfrac{-(b-2 \beta)}{c-2 \gamma}+\cfrac{(a- 2 \alpha):(c- 2 \gamma)}{\cfrac{-(b-3 \beta)}{c- 3 \gamma}+\cfrac{(a- 3 \alpha):(c- 3 \gamma)}{\cfrac{-(b-4 \beta)}{c-\gamma}+\text{etc.}}}}
\]
or having removed the partial fractions one finds

\[
\frac{(c-\gamma) \int PRdx}{\int Pdx}=-(b-\beta)+\cfrac{(a-\alpha)(c- 2 \gamma)}{-(b-2 \beta)+\cfrac{(a-2 \alpha)(c-3 \gamma)}{-(b-3 \beta)+\cfrac{(a-3\alpha)(c- 4 \gamma)}{-(b-4 \beta)+\text{etc.}}}}
\]
So one has two continued fractions, whose value is the same, namely $\frac{\int PRdx}{\int Pdx}$.

\paragraph*{§52}

But it is especially important here to define suitable functions $P$ and $R$ of $x$, that it it

\[
(a+n \alpha) \int PR^ndx =(b+n \beta) \int PR^{n+1}dx+(c+n \gamma) \int P R^{n+2}dx,
\]
at least in that case, in which after each integration one puts $x=1$. Now let us put that  it is in general

\[
(a+n \alpha) \int PR^ndx +R^{n+1}S =(b+n \beta) \int PR^{n+1}dx+(c+n \gamma) \int P R^{n+2}dx
\]
and that $R^{n+1}S$ is a function of such a kind of $x$ that vanishes so for $x=0$ as for $x=1$. Having taken differentials and having divided by $R^n$ it will be

\[
(a+n \alpha Pdx+RdS+(n+1)SdR=(b+n \beta)PRdx+(c+n \gamma)PR^2dx;
\]
this equation, because it has to hold in every case, no matter what $n$ is, is resolved into these two equations

\[
aPdx+RdS+SdR=bPRdx+cPR^2dx
\]
and

\[
\alpha Pdx+SdR= \beta PRdx+\gamma PR^2dx.
\]
From these equations one finds 

\[
Pdx=\frac{RdS+SdR}{bR+cR^2-a}=\frac{SdR}{\beta R+\gamma R^2-\alpha},
\]
whence it is

\[
\frac{dS}{S}=\frac{(b-\beta)RdR+(c-\gamma)R^2dR-(a-\alpha)dR}{\beta R^2+\gamma R^3- \alpha R}
\]
\[
=\frac{(a-\alpha)dR}{\alpha R}+\frac{(\alpha b- \beta a)dR+(\alpha c- \gamma a)RdR}{\alpha(\beta R+\gamma R^2-\alpha)}.
\]
Hence $S$ is defined by $R$ from this equation; but having found $S$, it will be

\[
P=\frac{SdR}{(\beta R- \gamma R^2- \alpha)dx}
\]
and hence the formulas $\int Pdx$ and $\int PRdx$ will be found, and these in turn determine the value of continued fraction given above.

\paragraph*{§53}
 Therefore, because the quantity $R$ is  not defined by $x$, one can take any function of $x$ for it. But because the nature of the question requires  $R^{n+1}S$ to vanish so for $x=0$ as for $x=1$, these conditions actually define the properties of the function of $R$. Further, one has to take into account that the integrals $\int PR^n dx$, having put $x=1$ after the integration, should obtain a finite value; hence if these integrals  would become either  $0$ or $\infty$ in this case, then it would be difficult to calculate the value of $\frac{\int PRdx}{\int Pdx}$. The first inconvenience is  avoided by assuming a function of such a kind for $R$ that $PR^n$ is a strictly positive function for $0<x<1$. But it is more difficult to achieve that $\int PR^ndx$ does not become infinite for $x=1$.  It will indeed be convenient,  to separate the cases, in which $n$ is a positive number from those, where $n$ is a negative number, because very often, if these conditions are fulfilled, while $n$ is a positive number,  they are not fulfilled in the remaining cases at the same time. But if the prescribed conditions are only fulfilled in the cases, in which $n$ is a positive number, then only the value of the first continued fraction can be exhibited; but the value of the second fraction on the other hand can only be given, if  the conditions were fulfilled, while $n$ is a negative number.
 
\paragraph*{§54}

Let us start the explanation of this method to find the values of continued fractions from the example treated before, and hence  first let this continued fraction be propounded 

\[
r+\cfrac{fh}{r+\cfrac{(f+r)(h+r)}{r+\cfrac{(f+2r)(h+2r)}{r+\text{etc.}}}}
\]
whose value was seen § $34$ to be

\[
\frac{h(f-r) \int y^{h+r-1}dy:\sqrt{1-y^{2r}}-f(h-r)\int y^{f+r-1}dy:\sqrt{1-y^{2r}}}{f\int y^{f+r-1}dy:\sqrt{1-y^{2r}}-h\int y^{h+r-1}dy:\sqrt{1-y^{2r}}}.
\]
Therefore, compare this continued fraction to this general one

\[
\frac{a\int Pdx}{\int PRdx}=b+\cfrac{(a+\alpha)c}{b+\beta+\cfrac{(a+2 \alpha)(c+\gamma)}{b+2 \beta+\cfrac{(a+3 \alpha)(c+2 \gamma)}{b+3 \beta+\text{etc.}}}}
\]
and it will be

\[
b=r, \quad \beta =0, \quad \alpha =r, \quad \gamma=r, \quad a=f-r, \quad c=h.
\]
Having substituted this value it will result

\[
\frac{ds}{S}= \frac{rRdR+(h-r)R^2dR-(f-2r)dR}{rR^3-rR}=\frac{(f-2r)dR}{rR}+\frac{rdR+(h-f+r)RdR}{r(R^2-1)}
\]
and by integrating

\[
\log S= \frac{f-2r}{r}\log R+\frac{h-f}{2r}\log (R+1)+\frac{h-f+2r}{2r}\log (R-1)+\log C
\]
or

\[
S=CR^{\frac{f-2r}{r}}(R^2-1)^{\frac{h-f}{2r}}(R-1).
\]
Therefore, it will hence be

\[
R^{n+1}S=CR^{\frac{f-(n-1)r}{r}}(R^2-1)^{\frac{h-f}{2r}}(R-1)
\]
and

\[
Pdx= \frac{CR^{\frac{f-2r}{r}}(R^2-1)^{\frac{h-f}{2r}}dR}{r(R+1)}.
\]

\paragraph*{§55}

But because $R^{n+1}S$ has to vanish in two cases, so for $x=0$ and $x=1$, and that, no matter which positive number is substituted for $n$ (for, to take into account  negative values is not necessary), let us put that $f$, $h$ and $r$ are positive numbers and $h>f$; this can certainly be assumed without loss of generality, if was not $f=h$; further, let it be $f>r$. Having constituted all this, it is manifest that the formula $R^{n+1}S$ vanishes in two cases, of course if $R=0$ and if $R=1$; and this  is also true, if it is $f=h$. If it is $f>r$, one can put $R=x$ and it will be

\[
Pdx=\frac{x^{\frac{f-2r}{r}}(1-x^2)^{\frac{h-f}{2r}}dx}{1+x}
\]
after having determined the constant $C$. Hence the value of the propounded continued fraction will be

\[
=(f-r)\frac{\int\frac{x^{\frac{f-2r}{r}}(1-x^2)^{\frac{h-f}{2r}}dx}{1+x}}{\int\frac{x^{\frac{f-r}{r}}(1-x^2)^{\frac{h-f}{2r}}dx}{1+x}}
\]
But for $x=y^r$ value in question will be

\[
= \frac{(f-r)\int y^{f-r-1}(1-y^{2r})^{\frac{h-f}{2r}}dy:(1+y^r)}{\int y^{f-1}(1-y^{2r})^{\frac{h-f}{2r}}dy:(1+y^r)}.
\]

\paragraph*{§56}

Hence we obtained another expression containing the value of this continued fraction

\[
r+\cfrac{fh}{r+\cfrac{(f+r)(h+r)}{r+\text{etc.}}},
\]
which expression, even though it contains integral formulas, nevertheless differs from the expression found before. Hence this last expression is only valid, if it is $f>r$; but for $h$ one has to take the smaller of the two quantities $f$ and $h$, if they were not equal, of course. But nevertheless, if also $f$ was was smaller than $r$, the value of the continued fraction can be exhibited by considering this one

\[
r+\cfrac{(f+r)(h+r)}{r+\cfrac{(f+2r)(h+2r)}{r+\text{etc.}}},
\]
whose value will be

\[
=\frac{f\int y^{f-1}(1-y^{2r})^{\frac{h-f}{2r}}dy:(1+y^r)}{\int y^{f+r-1}(1-y^{2r})^{\frac{h-f}{2r}}dy:(1+y^r)},
\]
which does not require any restrictions. Hence  having put this value $=V$ the value of the propounded continued fraction will be $=r+\frac{fh}{V}$.

\paragraph*{§57}

That case, in which $f=h$, which had been found before by a peculiar method and its value was detected in § $34$ to be

\[
=\frac{1-(h-r)\int x^{h-1}dx:(1+x^r)}{\int x^{h-1}dx:(1+x^r)}=\frac{(h-r)\int x^{h-r-1}dx:(1+x^r)}{\int x^{h-1}dx:(1+x^r)},
\]
follows from this last expression immediately; hence for $f=h$ the expression found in § $55$ will become this one

\[
\frac{(h-r)\int y^{h-r-1}dy:(1+y^r)}{\int y^{h-1}dy:(1+y^r)};
\]
this is completely the same expression, whence the agreement of the both general expressions is seen clearly. Here it is certainly possible to assume $h>r$, because the cases, in which this is not true, are very easily reduced to those in which this condition is met, as it was just shown.

\paragraph*{§58}

In order to understand this agreement of the two expression in every case, we have to mention the following lemma in advance, which was already proved by others: If we have this series
\[
1+\frac{p}{q+s}+\frac{p(p+s)}{(q+s)(q+2s)}+\frac{p(p+s)(p+2s)}{(q+s)(q+2s)(q+3s)}+\text{etc.},
\]
in which the quantities $p$, $q$ and $s$ are positive and $q>p$, the sum of this infinite series  will be $=\frac{q}{q-p}$. But this lemma can be proved by my general method to sum series in the following way. Just consider this series

\[
x^q+\frac{p}{q+s}x^{q+s}+\frac{p(p+s)}{(q+s)(q+2s)}x^{q+2s}+\text{etc.};
\]
call its sum $z$, and  by differentiating it will be

\[
\frac{dz}{dx}=qx^{q-1}+px^{q+s-1}+\frac{p(p+s)}{q+s}x^{q+2s-1}+\text{etc.}
\]
and

\[
x^{p-q-z}dz=qx^{p-s-1}dx+px^{p-1}dx+\frac{p(p+s)}{q+s}x^{p+s-1}dx+\text{etc.};
\]
integrating this equation gives

\[
\int x^{p-q-z}dz= \frac{qx^{p-s}}{p-s}+x^p+\frac{px^{p+s}}{p+s}+\text{etc.}=\frac{qx^{p-s}}{p-s}+x^{p-q}z.
\]
Differentiating this equation on the other hand gives 

\[
x^{p-q-s}ds=qx^{p-s-1}dx+x^{p-q}dx+(p-q)x^{p-q-1}zdx
\]
or

\[
dz(1-x^s)+(q-p)x^{s-1}zdx=qx^{q-1}dx
\]
or

\[
dz+\frac{(q-p)x^{s-1}zdx}{1-x^s}=\frac{qx^{q-1}dx}{1-x^s},
\]
whose integral is

\[
\frac{z}{(1-x^s)^{\frac{q-p}{s}}}=q\int \frac{x^{p-1}dx}{(1-x^s)^{\frac{q-p+s}{s}}}=\frac{qx^q}{(q-p)(1-x^s)^{\frac{q-p}{s}}}-\frac{pq}{q-p}\int \frac{x^{q-1}dx}{(1-x^s)^{\frac{q-p}{s}}},
\]
whence it will be

\[
z=\frac{qx^q}{q-p}-\frac{pq(1-x^s)^{\frac{q-p}{s}}}{q-p}\int \frac{x^{q-1}dx}{(1-x^s)^{\frac{q-p}{s}}}.
\]
Therefore, for $x=1$ it will be

\[
z=\frac{q}{q-p}=1+\frac{p}{q+s}+\frac{p(p+s)}{(q+s)(q+2s)}+\text{etc.},
\]
which is the demonstration of the stated lemma, whence at the same  time it is understood that it is only correct, if $q>p$.

\paragraph*{§59}

Therefore,  we expressed the value of the continued fraction

\[
r+\cfrac{fh}{r+\cfrac{(f+r)(h+r)}{r+\cfrac{(f+2r)(h+2r)}{r+\text{etc.}}}}
\]
 in two ways; first

\[
=\frac{h(f-r)\int y^{h+r-1}dy:\sqrt{1-y^{2r}}-f(h-r)\int y^{f+r-1}dy:\sqrt{1-y^{2r}}}{f \int y^{f+r-1}dy: \sqrt{1-y^{2r}}-h\int y^{h+r-1}dy:\sqrt{1-y^{2r}}},
\]
and second, as it was found in § $56$,

\[
=r+\frac{h \int y^{f+r-1}dy(1-y^{2r})^{\frac{h-f}{2r}}:(1+y^r)}{\int y^{f-1} dy(1-y^{2r})^{\frac{h-f}{2r}}:(1+y^r)}.
\]
It will be worth one's while to demonstrate the agreement of these expressions. Therefore, because it is

\[
\frac{1}{1+y^r}=\frac{1-y^r}{1-y^{2r}},
\]
it will be

\[
\int y^{f-1}dy(1-y^{2r})^{\frac{h-f}{2r}}:(1+y^r)=\int y^{f-1}dy(1-y^{2r})^{\frac{h-f-2r}{2r}}-\int y^{f+r-1}dy(1-y^{2r})^{\frac{h-f-2r}{2r}}
\]
and

\[
\int y^{f+r-1}dy(1-y^{2r})^{\frac{h-f}{2r}}:(1+y^r)
\]
\[
=\int y^{f+r-1}dy(1-y^{2r})^{\frac{h-f-2r}{2r}}-\int y^{f+2r-1}dy(1-y^{2r})^{\frac{h-f-2r}{2r}}
\]
\[
=\int y^{f+r-1}dy(1-y^{2r})^{\frac{h-f-2r}{2r}}-\frac{f}{h}\int y^{f-1}dy(1-y^{2r})^{\frac{h-f-2r}{2r}}.
\]
Now just put

\[
\frac{\int y^{f+r-1}dy(1-y^{2r})^{\frac{h-f-2r}{2r}}}{\int y^{f-1}dy(1-y^{2r})^{\frac{h-f-2r}{2r}}}=V;
\]
the last value of the continued fraction will be

\[
=r+\frac{hV-f}{1-V}.
\]
Now furthermore put

\[
\frac{\int y^{h+r-1}dy:\sqrt{1-y^{2r}}}{\int y^{f+r-1}dy :\sqrt{1-y^{2r}}}=W;
\]
the first value will be

\[
=\frac{h(f-r)W-f(h-r)}{f-hW};
\]
and from the equality of the two values it follows that it will be

\[
V=\frac{f}{hW},
\]
so that it is

\[
\frac{\int y^{f+r-1}dy(1-y^{2r})^{\frac{h-f-2r}{2r}}}{\int y^{f-1}dy(1-y^{2r})^{\frac{h-f-2r}{2r}}}=\frac{f\int y^{f+r-1}dy:\sqrt{1-y^{2r}}}{h\int y^{h+r-1}dy:\sqrt{1-y^{2r}}};
\]
the reason for this equality is known from the theorems exhibited in the preceding dissertation\footnote{Euler means E122 again.}; hence  by one of those theorems it will be

\[
\frac{\int y^{f+r-1}dy(1-y^{2r})^{\frac{h-f-2r}{2r}}}{\int y^{f+2r-1}dy(1-y^{2r})^{\frac{h-f-2r}{2r}}}=\frac{\int y^{f+r-1}dy:\sqrt{1-y^{2r}}}{\int y^{h+r-1}dy:\sqrt{1-y^{2r}}}.
\]

\paragraph*{§60}

Now let us consider this continued fraction

\[
2r+\cfrac{fh}{2r+\cfrac{(f+r)(h+r)}{2r+\cfrac{(f+2r)(h+2r)}{2r+\text{etc.}}}}
\]
whose value was found above in § $35$ to be

\[
=\frac{2(f-r)(h-r)\int y^{f-1}dy:\sqrt{1-y^{2r}}-h(f+h-3r)\int y^{h+r-1}dy:\sqrt{1-y^{2r}}}{2h \int y^{h+r-1}dy:\sqrt{1-y^{2r}}-(f+h-r)\int y^{f-1}dy:\sqrt{1-y^{2r}}}
\]
If now that continued fraction is compared to this one

\[
\frac{a \int Pdx}{\int PRdx}=b+\cfrac{(a+\alpha)c}{b+\beta+\cfrac{(a+2 \alpha)(c+\gamma)}{b+2 \beta+\cfrac{(a+3 \alpha)(c+2 \gamma)}{b+2 \beta+\text{etc.}}}}
\]
it will be

\[
b=2r, \quad \beta=0, \quad \alpha =r, \quad \gamma=r, \quad a=f-r \quad \text{and} \quad c=h.
\]
Hence from § $52$ one will have

\[
\frac{dS}{S}=\frac{(f-2r)dR}{rR}+\frac{2rdR+(h-f+r)RdR}{r(R^2-1)}
\]
and by integrating

\[
S=CR^{\frac{f-2r}{2r}}(R^2-1)^{\frac{h-f-r}{2r}}(R-1)^2,
\]
whence it is

\[
Pdx= \frac{C}{r}R^{\frac{f-2r}{r}}(R^2-1)^{\frac{h-f-2r}{2r}}(R-1)^2dR
\]
and

\[
R^{n+1}S=CR^{\frac{f+(n-1)r}{r}}(R^2-1)^{\frac{h-f-r}{2r}}(R-1)^2,
\]
which expression vanishes in two cases, for $R=0$ and for $R=1$, if only $f>r$ and $h-3r>f$, which conditions can always be fulfilled.

\paragraph*{§61}
Now let $R=x$ and having determined the constant $C$ it will be
 
\[
Pdx=x^{\frac{f-2r}{r}}dx(1-x^2)^{\frac{h-f-3r}{2r}}(1-x)^2
\]
or having put $R=x=y^r$ it will be

\[
Pdx=y^{f-r-1}dy(1-y^{2r})^{\frac{h-f-3r}{2r}}(1-y^r)^2;
\]
using this  value of the propounded continued fraction it will be

\[
\frac{a\int Pdx}{\int PRdx}=\frac{(f-r)\int y^{f-r-1}dy(1-y^{2r})^{\frac{h-f-3r}{2r}}(1-y^r)^2}{\int y^{f-1}dy(1-y^{2r})^{\frac{h-f-3r}{2r}}(1-y^r)^2};
\]
this expression, applying the theorems of the preceding dissertation\footnote{E122},   will be reduced to the  first form by expanding the square $(1-y^r)^2$; having done this  both integral formulas will be resolved into two simpler ones. But I will show the reduction in the following more general example.

\paragraph*{§62}

If one has this integral formula

\[
\int y^{m-1}dy(1-y^{2r})^{\varkappa}(1-y^r)^n
\] 
and $(1-y^r)^n$ is expanded into the series

\[
1-ny^r+\frac{n(n-1)}{1 \cdot 2}y^{2r}-\text{etc.},
\]
by taking every second term the propounded integral formula will be reduced to the two following ones 
\[
\int y^{m-1}dy(1-y^{2r})^{\varkappa}\bigg(1+\frac{n(n-1)}{1 \cdot 2}\cdot \frac{m}{p}+\frac{n(n-1)(n-2)(n-3)}{1 \cdot 2 \cdot 3 \cdot 4}\cdot \frac{m(m+2r)}{p(p+2r)}+\text{etc.}\bigg)
\]
\[
-\int y^{m+r-1}dy(1-y^{2r})^{\varkappa} \left\{
\begin{aligned}
& \quad \quad \quad \quad \quad \quad \quad n+\frac{n(n-1)(n-2)}{1 \cdot 2 \cdot 3}\cdot \frac{m+r}{p+r}\\
& +\frac{n(n-1)(n-2)(n-3)(n-4)}{1 \cdot 2 \cdot 3\cdot 4 \cdot 5} \cdot \frac{(m+r)(m+3r)}{(p+r)(p+3r)}+\text{etc.}
\end{aligned}
\right\}
\]
after having put

\[
m+2 \varkappa r+2r=p
\]
for the sake of brevity. Therefore, if as in the preceding case it was $n=2$, it will be

\[
\int y^{m-1}dy(1-y^{2r})^{\varkappa}(1-y^r)^2= \frac{m+p}{p}\int y^{m-1}dy(1-y^{2r})^{\varkappa}-2\int y^{m+r-1}dy(1-y^{2r})^{\varkappa}.
\]
From this one will have

\[
\frac{a\int Pdx}{\int PRdx}=\frac{\frac{(f-r)(f+h-3r)}{h-2r}\int y^{f-r-1}dy(1-y^{2r})^{\frac{h-f-3r}{2r}}-2(f-r)\int y^{f-1}dy(1-y^{2r})^{\frac{h-f-3r}{2r}}}{\frac{f+h-r}{h-r}\int y^{f-1}dy(1-y^{2r})^{\frac{h-f-3r}{2r}}-2\int y^{f+r-1}dy(1-y^{2r})^{\frac{h-f-3r}{2r}}}
\]
\[
=\frac{h(f+h-3r)\int y^{f-r-1}dy(1-y^{2r})^{\frac{h-f-r}{2r}}-2(f-r)(h-r)\int y^{f-1}dy(1-y^{2r})^{\frac{h-f-r}{2r}}}{(f+h-r)\int y^{f-1}dy(1-y^{2r})^{\frac{h-f-r}{2r}}-2h\int y^{f+r-1}dy(1-y^{2r})^{\frac{h-f-r}{2r}}},
\]
which expression, because it has to be equal to the one found above in § $35$, will yield this equation

\[
\frac{\int y^{f+r-1}dy(1-y^{2r})^{\frac{h-f-r}{2r}}}{\int y^{f-1}dy(1-y^{2r})^{\frac{h-f-r}{2r}}}=\frac{\int y^{h+r-1}dy:\sqrt{1-y^{2r}}}{\int y^{f-1}dy:\sqrt{1-y^{2r}}};
\]
the reason for this is indeed contained in the theorems of the preceding dissertation\footnote{Euler refers to E122 again.}.

\paragraph*{§63}

Now let us vice versa take given functions for $P$ and $R$ and let us form continued fractions using them and let us put

\[
P=x^{m-1}(1-x^r)^n(p+qx^r)^{\varkappa} \quad \text{and} \quad R=x^r.
\]
But because it has to be

\[
(a+\nu \alpha)\int PR^{\nu}dx=(b+\nu \beta)\int PR^{\nu +1}dx+(c+\nu \gamma)\int PR^{\nu+2}dx
\]
and hence because of the given $P$ and $R$  from § $52$ it is

\[
S=\frac{1}{r}x^{m-r}(1-x^r)^n(p+qx^r)^{\varkappa}(\gamma x^{2r}+\beta x^r-\alpha),
\]
it will be

\[
\frac{dS}{S}=\frac{(m-r)dx}{x}+\frac{nrx^{r-1}dx}{-1+x^r}+\frac{\varkappa q r x^{r-1}dx}{p+qx^r}+\frac{2 \gamma rx^{2r-1}dx+\beta r x^{r-1}dx}{\gamma x^{2r}+\beta x^r-\alpha}
\]
\[
=\frac{(a-\alpha)rdx}{\alpha x}+\frac{(\alpha b-\beta a)rx^{r-1}dx+(\alpha c-\gamma a)rx^{2r-1}dx}{\alpha(\gamma x^{2r}+\beta x^r-\alpha)}.
\]
Now let be

\[
(p+qx^r)(x^r-1)=\gamma x^{2r}+\beta x^r-\alpha;
\]
it will be

\[
\gamma =q, \quad \beta =p-q \quad \text{and} \quad \alpha =p.
\]
Furthermore, let it be

\[
\frac{(a-\alpha)r}{\alpha}=m-r;
\]
it will be

\[
a=\frac{mp}{r}.
\]
Hence further it will be

\[
nqr+\varkappa qr+2 qr=\frac{cpr-mpq}{p}
\]
or

\[
c=\frac{mq}{r}+nq+(\varkappa +2)q
\]
and finally

\[
b= \frac{m(p-q)}{r}+(n+1)p-(\varkappa +1)q.
\]
So as long as $m$ and $n+1$ were positive numbers that $R^{\nu+1}S$ vanishes so for $x=$ as for $x=1$, the following expression will result

\[
\frac{\int x^{m+r-1}dx(1-x^r)^n(p+qx^r)^{\varkappa}}{\int x^{m-1}dx(1-x^r)^n(p+qx^r)^{\varkappa}}=\frac{\int PRdx}{\int Pdx},
\]
which will therefore be equal to this continued fraction

\[
\cfrac{mp}{m(p-q)+(n+1)pr-(\varkappa+1)qr+\cfrac{pq(m+r)(m+nr+(\varkappa+2)r)}{m(p-q)+(n+2)pr-(\varkappa +2)qr+}}
\]
\[
\cfrac{pq(m+2r)(m+(n+1)r+(\varkappa+2)r)}{m(p-q)+(n+3)pr-(\varkappa +3)qr+\text{etc.}}
\]

\paragraph*{§64}

To simplify this continued fraction just put

\[
m+nr+r=a, \quad m+\varkappa r+r=b \quad \text{and} \quad m+nr+\varkappa r+r=c;
\]
it will be

\[
\varkappa=\frac{c-a}{r}, \quad n=\frac{c-b}{r} \quad \text{and} ßquad m=a+b-c-r
\]
and therefore it will be

\[
\cfrac{p(a+b-c-r)}{ap-bq+\cfrac{pq(a+b-c)(c+r)}{(a+r)p-(b+r)q+\cfrac{pq(a+b-c+r)(c+2r)}{(a+2r)p-(b+2r)q+\cfrac{pq(a+b-c+2r)(c+3r)}{(a+3r)p-(b+3r)q+\text{etc.}}}}}
\]
\[
=\frac{\int x^{a+b-c-1}dx(1-x^r)^{\frac{c-b}{r}}(p+qx^r)^{\frac{c-a}{r}}}{\int x^{a+b-c-r-1}dx(1-x^r)^{\frac{c-b}{r}}(p+qx^r)^{\frac{c-a}{r}}}
\]
having put $x=1$ after each integration. But it is required that

\[
a+b-c-r \quad \text{and} \quad c-b+r
\]
are positive numbers. But if for the sake of brevity it is put

\[
a+b-c-r=g,
\] 
it will be

\[
\frac{\int x^{g-r-1}dx(1-x^r)^{\frac{c-b}{r}}(p+qx^r)^{\frac{c-a}{r}}}{\int x^{g-1}dx(1-x^r)^{\frac{c-b}{r}}(p+qx^r)^{\frac{c-a}{r}}}
\]
\[
=\frac{pg}{ap-bq+\cfrac{pq(c+r)(g+r)}{(a+r)p-(b+r)q+\cfrac{pq(c+2r)(g+2r)}{(a+2r)p-(b+2r)q+\text{etc.}}}}
\]
which equation extends very far and contains all continued fractions found up to now.

\paragraph*{§65}

If the quantities $c$ and $g$ are permuted, the following continued fraction will result

\[
\frac{pc}{ap-bq+\cfrac{pq(c+r)(g+r)}{(a+r)p-(b+r)q+\cfrac{pq(c+2r)(g+2r)}{(a+2r)p-(b+2r)q+\text{etc.}}}},
\]
whose value will therefore be

\[
\frac{\int x^{c+r-1}dx(1-x^r)^{\frac{g-b}{r}}(p+qx^r)^{\frac{g-a}{r}}}{\int x^{c-1}dx(1-x^r)^{\frac{g-b}{r}}(p+qx^r)^{\frac{g-a}{r}}}.
\]
Because these continued fractions have a given ratio, of course $g$ to $c$, the following theorem will follow from this after having resubstituted  the value of $g$

\[
\frac{c\int x^{a+b-c-1}dx(1-x^r)^{\frac{c-b}{r}}(p+qx^r)^{\frac{c-a}{r}}}{\int x^{a+b-c-r-1}dx(1-x^r)^{\frac{c-b}{r}}(p+qx^r)^{\frac{c-a}{r}}}
\]
\[
=\frac{(a+b-c-r)\int x^{c+r-1}dx(1-x^r)^{\frac{a-c-r}{r}}(p+qx^r)^{\frac{b-c-r}{r}}}{\int x^{c-1}dx(1-x^r)^{\frac{a-c-r}{r}}(p+qx^r)^{\frac{b-c-r}{r}}}.
\]
This very general formula contains many extraordinary particular reductions. So  for the sake of an example let it be $b=c+r$; it will be

\[
\frac{a\int x^{a+r-1}dx(p+qx^r)^{\frac{c-a}{r}}:(1-x^r)}{\int x^{a-1}dx(p+qx^r)^{\frac{c-a}{r}}:(1-x^r)}=\frac{a\int x^{a+r-1}dx(1-x^r)^{\frac{a-c-r}{r}}}{\int x^{c-1}dx(1-x^r)^{\frac{c-a}{r}}}=c,
\]
whence it follows that it will be

\[
\int \frac{x^{a+r-1}dx(p+qx^r)^{\frac{c-a}{r}}}{1-x^r}=\int \frac{x^{a-1}dx(p+qx^r)^{\frac{c-a}{r}}}{1-x^r}.
\]
One will therefore have this more general extending theorem

\[
\int \frac{x^{m-1}dx(p+qx^r)^{\varkappa}}{1-x^r}=\int \frac{x^{n-1}dx(p+qx^r)^{\varkappa}}{1-x^r},
\]
where, after having integrated in such a way that the integrals vanish for $x=0$, $x$ is understood to become $=1$. And indeed only that case, in which $q+p=0$ and in which the inconvenience occurs, is excluded.

\paragraph*{§66}

The continued fractions we found up to now by interpolation, reduce to this that the partial denominators all have the same value. To express those continued fractions by this  general form, put $p=q=1$ and this continued fraction will result

\[
\cfrac{cg}{a-b+\cfrac{(c+r)(g+r)}{a-b+\cfrac{(c+2r)(g+2r)}{a-b+\cfrac{(c+3r)(g+3r)}{a-b+\text{etc.}}}}}= \frac{c\int x^{g+r-1}dx(1-x^r)^{\frac{c-b}{r}}(1+x^r)^{\frac{c-a}{r}}}{\int x^{g-1}dx(1-x^r)^{\frac{c-b}{r}}(1+x^r)^{\frac{c-a}{r}}}
\]
or the value of the same will also be

\[
=\frac{g\int x^{c+r-1}dx(1-x^r)^{\frac{g-b}{r}}(1+x^r)^{\frac{g-a}{r}}}{\int x^{c-1}dx(1-x^r)^{\frac{g-b}{r}}(1+x^r)^{\frac{g-a}{r}}}
\]
while it is $g=a+b-c-r$. Just put

\[
a-b=s;
\] 
because of

\[
a+b=c+g+r
\]
it will be

\[
a=\frac{c+g+r+s}{2} \quad \text{and} \quad b=\frac{c+g+r-s}{2},
\]
whence it will be

\[
\frac{cg}{s+\cfrac{(c+r)(g+r)}{s+\cfrac{(c+2r)(g+2r)}{s+\text{etc.}}}}
\]
\[
=\frac{c \int x^{g+r-1}dx(1-x^{2r})^{\frac{c-g-r-s}{2r}}(1-x^r)^{\frac{s}{r}}}{ \int x^{g-1}dx(1-x^{2r})^{\frac{c-g-r-s}{2r}}(1-x^r)^{\frac{s}{r}}}=\frac{g \int x^{c+r-1}dx(1-x^{2r})^{\frac{g-c-r-s}{2r}}(1-x^r)^{\frac{s}{r}}}{c \int x^{c-1}dx(1-x^{2r})^{\frac{g-c-r-s}{2r}}(1-x^r)^{\frac{s}{r}}}.
\]

\paragraph*{§67}

In order to get to the from in § $47$, let us put $2s$ instead of $s$ and let be $c=q$ and $g=r-q$; one will have this continued fraction

\[
\frac{q(r-q)}{2s+\cfrac{(q+r)(2r-q)}{2s+\cfrac{(q+2s)(3r-q)}{2s+\text{etc.}}}},
\]
whose value will therefore be either

\[
=\frac{q \int x^{2r-q-1}dx(1-x^{2r})^{\frac{q-r-s}{r}}(1-x^r)^{\frac{2s}{r}}}{\int x^{r-q-1}dx(1-x^{2r})^{\frac{q-r-s}{r}}(1-x^r)^{\frac{2s}{r}}}
\]
or

\[
=\frac{(r-q)\int x^{q+r-1}dx(1-x^{2r})^{\frac{-q-s}{r}}(1-x^r)^{\frac{2s}{r}}}{\int x^{q-1}dx(1-x^{2r})^{\frac{-q-s}{r}}(1-x^r)^{\frac{2s}{r}}}.
\]
The value of the same continued fraction was found before to be

\[
=\frac{(q+s)\int y^{r+s+q-1}dy:\sqrt{1-y^{2r}}}{\int y^{r+s-q-1}dy:\sqrt{1-y^{2r}}}-s.
\]
Therefore, these integral formulas will be equal to each other; this is an outstanding theorem.

\paragraph*{§68}

Let us, as we did in § $48$, set $r=2$ and $q=1$; it will be

\[
\frac{(1+s)\int y^{s+2}dy:\sqrt{1-y^4}}{\int y^sdy:\sqrt{1-y^4}}-s=\frac{\int x^2dx(1-x^4)^{\frac{-s-1}{2}}(1-x^2)^s}{\int dx(1-x^4)^{-\frac{-s-1}{2}}(1-x^2)^s},
\]
which equality is obvious, if $s=0$; and even in the cases, in which $s$ is an odd integer, the equality is easily shown. So if it is $s=1$, the last formula will be

\[
\frac{\int xxdx:(1+xx)}{\int dx:(1+xx)}=\frac{x-\int dx:(1+xx)}{\int dx:(1+xx)}=\frac{4-\pi}{\pi}
\]
for $x=1$. The first formula on the other hand will give

\[
\frac{2\int y^3dy:\sqrt{1-y^4}}{\int ydy:\sqrt{1-y^4}}-1=\frac{4}{\pi}-1=\frac{4-\pi}{\pi};
\]
hence it agrees completely with the preceding. But if $s$ is an even number, the agreement of the both expressions is easily seen by expansion of $(1-xx)^s$.

\paragraph*{§69}

But except for the continued fractions found up to now the  general formula contains innumerable others; and it will be helpful to have expanded some of them. Therefore, let it be $g=c$ and the value of this continued fraction

\[
\cfrac{c^2}{s+\cfrac{(c+r)^2}{s+\cfrac{(c+2r)^2}{s+\text{etc.}}}}
\] 
will be

\[
\frac{c\int x^{c+r-1}dx(1-x^r)^{\frac{s}{r}}:(1-x^{2r})^{\frac{r+s}{2r}}}{\int x^{c-1}dx(1-x^r)^{\frac{s}{r}}:(1-x^{2r})^{\frac{r+s}{2r}}}.
\]
Now put $c=1$ and $r=1$ and it will be

\[
\cfrac{1}{s+\cfrac{4}{s+\cfrac{9}{s+\cfrac{16}{s+\text{etc.}}}}}=\frac{\int xdx(1-x)^s:(1-xx)^{\frac{s+1}{2}}}{\int dx(1-x)^s:(1-xx)^{\frac{s+1}{2}}};
\]
let us investigate the values of this expression for various values of $s$. Therefore, having put the value of this expression $=V$, it will be as follows:

\begin{alignat*}{9}
& \text{if} \quad &&s=0, \quad && \\
& &&V=\frac{\int xdx:\sqrt{1-xx}}{\int dx:\sqrt{1-xx}}&&=\frac{1}{2 \int dy:(1+yy)};\\
& \text{if} \quad &&s=2, \quad &&  \\
& &&V=\frac{2\int dx:\sqrt{1-xx}-3 \int xdx:\sqrt{1-xx}}{2\int xdx:\sqrt{1-xx}-\int dx:\sqrt{1-xx}}&&=\frac{1}{2 \int y^2dy:(1+yy)}-2;\\
& \text{if} \quad &&s=4, \quad &&  \\
& &&V=\frac{19\int xdx:\sqrt{1-xx}-12 \int dx:\sqrt{1-xx}}{3\int dx:\sqrt{1-xx}-4\int xdx:\sqrt{1-xx}}&&=\frac{1}{2 \int y^4dy:(1+yy)}-4.
\end{alignat*}
But in general it will be

\[
V=\frac{1}{2 \int y^sdy:(1+yy)}-s;
\]
from this formula it is obvious, if $S$ was an even integer that the integral can be reduced to the quadrature of the circle; but, if $s$ was odd, the integrals can be expressed by logarithms.

\paragraph*{§70}

Now let this continued fraction be propounded 

\[
1+\cfrac{1}{2+\cfrac{4}{3+\cfrac{9}{4+\cfrac{16}{5+\cfrac{25}{6+\text{etc.}}}}}}
\]
Now compare this one to the form exhibited in § $64$ and it will be

\[
pqcg=1,
\]
\[
pq(c+r)(g+r)=4,
\]
\[
pq(c+2r)(g+2r)=9,
\]
\[
ap-bq=2 \quad \text{and} \quad (p-q)r=1,
\]
whence it will be

\[
p=\frac{\sqrt{5}+1}{2r}, \quad q= \frac{\sqrt{5}-1}{2r},
\]
\[
a= \frac{r(1+3\sqrt{5})}{2 \sqrt{5}} \quad \text{and} \quad b= \frac{r(3\sqrt{5}-1)}{2\sqrt{5}}; 
\]
after having substituted these values one will find the value of the propounded continued fraction to be

\[
=1+\frac{(\sqrt{5}-1)\int x^{2r-1}dx(1-x^r)^{\frac{1-\sqrt{5}}{2\sqrt{5}}}(1+\sqrt{5}+(\sqrt{5}-1)x^r)^{\frac{-\sqrt{5}-1}{2\sqrt{5}}}}{2\int x^{r-1}dx(1-x^r)^{\frac{1-\sqrt{5}}{2\sqrt{5}}}(1+\sqrt{5}+(\sqrt{5}-1)x^r)^{\frac{-\sqrt{5}-1}{2\sqrt{5}}}}.
\]
But, because of the irrational exponents, nothing special can be deduced from this formula.

\paragraph*{§71}

Whereas in these continued fractions the partial numerators are composed of two factors, I now proceed to continued fractions of such kind, in which these partial numerators constitute an arithmetic progression. Therefore, now, returning to § $50$, let it be $\gamma =0$ and $c=1$; it will be

\[
\frac{\int PRdx}{\int Pdx}=\cfrac{a}{b+\cfrac{a+\alpha}{b+\beta+\cfrac{a+2 \alpha}{b+ 2 \beta+\cfrac{a+3 \alpha}{b+3 \beta +\text{etc.}}}}}
\] 
But one has to take

\[
\frac{dS}{S}=\frac{(a-\alpha)dR}{\alpha R}+\frac{(\alpha b-\beta a)dR+\alpha RdR}{\alpha(\beta R-\alpha)}=\frac{(a-\alpha)dR}{\alpha R}+\frac{dR}{\beta}+\frac{(\alpha^2+\alpha \beta b-\beta ^2 a)dR}{\alpha \beta (\beta R-\alpha)},
\]
whence it is

\[
S=Ce^{\frac{R}{\beta}}R^{\frac{a-\alpha}{\alpha}}(\beta R-\alpha)^{\frac{\alpha ^2 +\alpha \beta b- \beta ^2 a}{\alpha \beta \beta}}.
\]
Put

\[
R=\frac{\alpha x}{\beta};
\]
it will be

\[
S=Ce^{\frac{\alpha x}{\beta \beta}}x^{\frac{a-\alpha}{\alpha}}(1-x)^{\frac{\alpha ^2+\alpha \beta b- \beta^2 a}{\alpha \beta \beta}}
\]
and $R^{n+1}S$ vanishes in two cases, of course so for $x=0$ as for $x=1$, as long

\[
\alpha ^2+\alpha \beta b >\beta ^2 a.
\]
Hence it will be

\[
Pdx=e^{\frac{\alpha x}{\beta \beta}}x^{\frac{a-\alpha}{\alpha}}dx(1-x)^{\frac{\alpha ^2+\alpha \beta b-\alpha \beta ^2- \beta^2 a}{\alpha \beta \beta}}
\]
and the value of the propounded continued fraction will be

\[
=\frac{\int PRdx}{\int Pdx}=\frac{\alpha \int e^{\frac{\alpha x}{\beta \beta}}x^{\frac{a}{\alpha}}dx(1-x)^{\frac{\alpha ^2+\alpha \beta b-\alpha \beta^2- \beta^2 a}{\alpha \beta \beta}}}{\beta \int e^{\frac{\alpha x}{\beta \beta}}x^{\frac{a-\alpha}{\alpha}}dx(1-x)^{\frac{\alpha ^2+\alpha \beta b-\alpha \beta^2- \beta^2 a}{\alpha \beta \beta}}}
\]
after having put $x=1$ after integration.

\paragraph*{§72}

To illustrate this case by an example, let

\[
a=1, \quad \alpha=1 \quad \text{and} \quad \beta=1;
\]
one will have this continued fraction

\[
\cfrac{1}{1+\cfrac{2}{2+\cfrac{3}{3+\cfrac{4}{4+\text{etc.}}}}}
\]
whose value will be

\[
\frac{\int e^xxdx}{\int e^xdx}=\frac{xe^x-e^x+1}{e^x-1}=\frac{1}{e-1}
\]
having put $x=1$. Hence it will be

\[
e=2+\cfrac{2}{2+\cfrac{3}{3+\cfrac{4}{4+\cfrac{5}{5+\text{etc.}}}}}
\]
by which expression one finds the decimal expansion of the number $e$, whose logarithm is $=1$, very quickly.

\paragraph*{§73}

Now let us put  $\beta=0$ in the  continued fraction given in § $71$, so that it 

\[
\frac{\int PRdx}{\int Pdx}=\cfrac{a}{b+\cfrac{a+\alpha}{b+\cfrac{a+2 \alpha}{b+\cfrac{a+3 \alpha}{b+\text{etc.}}}}}
\]
it will be

\[
\frac{dS}{S}=\frac{(a-\alpha)dR}{\alpha R}-\frac{bdR}{\alpha}-\frac{RdR}{\alpha}
\]
and hence

\[
S=CR^{\frac{a-\alpha}{\alpha}}e^{\frac{-2bR-RR}{2 \alpha}}.
\]
Now $R^{n+1}S$ vanishes in two cases; first, if $R=0$, then, if $R=\infty$, as long $a$ and $\alpha$ are positive numbers. Therefore, put

\[
R=\frac{x}{1-x}
\]
and it will be

\[
S=Cx^{\frac{a- \alpha}{\alpha}}:(1-x)^{\frac{a- \alpha}{\alpha}}e^{\frac{2bx-(2b-1)xx}{2 \alpha(1-x)^2}}.
\]
Because of

\[
dR=\frac{dx}{(1-x)^2}
\]
it will be

\[
\int Pdx= \int \frac{x^{\frac{a-\alpha}{\alpha}}dx}{(1-x)^{\frac{a+\alpha}{\alpha}}e^{\frac{2bx-(2b-1)xx}{2 \alpha(1-x)^2}}}
\]
and

\[
\int PRdx= \int \frac{x^{\frac{a}{\alpha}}dx}{(1-x)^{\frac{a+2\alpha}{\alpha}}e^{\frac{2bx-(2b-1)xx}{2 \alpha(1-x)^2}}}.
\]

\paragraph*{§74}

Finally,  in § $50$ let it be

\[
a=1, \quad c=1, \quad \alpha=0 \quad \gamma=0;
\]
it will be 

\[
\frac{\int PRdx}{\int Pdx}=\cfrac{1}{b+\cfrac{1}{b+\beta+\cfrac{1}{b+2 \beta+\cfrac{1}{b+3 \beta+\text{etc.}}}}}
\]
and

\[
\frac{dS}{S}=\frac{R^2dR+(b-\beta)RdR-dR}{\beta R^2},
\]
whence it will be

\[
S=e^{\frac{RR+1}{\beta}}R^{\frac{b-\beta}{\beta}} \quad \text{and} \quad Pdx=e^{\frac{RR+1}{\beta}}R^{\frac{b-2\beta}{\beta}}dR
\]
and

\[
PRdx=e^{\frac{RR+1}{\beta}}R^{\frac{b-\beta}{\beta}}dR.
\]
But $R$ has to be a function of $x$ of such a kind that $R^{n+1}$ vanishes so for $x=0$ as for $x=1$. But it is a lot more difficult to assign a function of such a kind  than for the remaining cases. Hence I will not try to resolve this case by the same method, but will solve it by another method now to be explained.

\paragraph*{§75}

I indeed mentioned this method to get to continued fractions already some time ago\footnote{In his first paper on continued fractions, E71.}, but because I only treated a particular case at that time, it will be convenient to explain it here in more detail. But this method does not use integral formulas as the preceding method, but relies on the resolution of a differential equation similar to that one considered by  \textsc{Riccati}. Of course, I consider this differential equation 

\[
ax^mdx+bx^{m-1}ydx+cy^2dx+dy=0,
\] 
which by putting

\[
x^{m+3}=0 \quad \text{and} \quad y=\frac{1}{cx}+\frac{1}{xxz} 
\]
goes over into this one

\[
\frac{-c}{m+3}t^{\frac{-m-4}{m+3}}dt-\frac{b}{m+3}t^{\frac{-1}{m+3}}zdt-\frac{ac+b}{(m+3)c}z^2dt+dz=0,
\]
which is similar to the first one. Hence if the value of $z$ was known in terms of $t$, at the same time $y$ will be known in terms of $x$. Therefore,  reduce this equation to another similar one in the same way by putting

\[
t^{\frac{2m+5}{m+3}}=u \quad \text{and} \quad z= \frac{-(m+3)c}{(ac+b)t}+\frac{1}{ttv}
\]
and just continue the reductions of this kind to infinity; having done this, if each value is substituted in  the preceding one, $y$ will be expressed in the following way

\[
y=Ax^{-1}+\cfrac{1}{-Bx^{-m-1}+\cfrac{1}{Cx^{-1}+\cfrac{1}{-Dx^{-m-1}+\cfrac{1}{Ex^{-1}+\cfrac{1}{-Fx^{-m-1}+\text{etc.}}}}}};
\]
the letters $A$, $B$, $C$, $D$ etc. on the other hand will obtain the following values

\begin{alignat*}{9}
&A&&=\frac{1}{c}, \\
&B&&=\frac{(m+3)c}{ac+b},\\
&C&&=\frac{(2m+5)(ac+b)}{c(ac-(m+2)b)},\\
&D&&=\frac{(3m+7)c(ac-(m+2)b)}{(ac+b)(ac+(m+3)b)},\\
&E&&=\frac{(4m+9)(ac+b)(ac+(m+3)b)}{c(ac-(m+2)b)(ac-(2m+4)b)},\\
&F&&=\frac{(5m+11)c(ac-(m+2)b)(ac-(2m+4)b)}{(ac+b)(ac+(m+3)b)(ac+(2m+5)b)}\\
& &&\text{etc.},
\end{alignat*}
which relations are more easily understood  by the following equations:

\begin{alignat*}{9}
&AB&&=\frac{m+3}{ac+b},\\
&BC&&=\frac{(m+3)(2m+5)}{ac-(m+2)b},\\
&CD&&=\frac{(2m+5)(3m+7)}{ac+(m+3)b},\\
&DE&&=\frac{(3m+7)(4m+9)}{ac-(2m+4)b},\\
&EF&&=\frac{(4m+9)(5m+11)}{ac+(2m+5)b},\\
&FG&&=\frac{(5m+11)(6m+13)}{ac-(3m+6)b}\\
& &&\text{etc.}
\end{alignat*}

\paragraph*{§76}

If now these values are substituted in the found continued fraction, one will then find

\[
cxy=1+\cfrac{(ac+b)x^{m+2}}{-(m+3)+\cfrac{(ac-(m+2)b)x^{m+2}}{(2m+5)+\cfrac{(ac+(m+3)b)x^{m+2}}{-(3m+7)+\cfrac{(ac-(2m+4)b)x^{m+2}}{(4m+9)+\text{etc.}}}}}
\]
From this expression it is clear that the propounded equation is absolutely integrable in the cases, in which $b$ becomes equal to a certain term of this progression

\[
-ac, \quad \frac{-ac}{m+3}, \quad \frac{-ac}{2m+5}, \quad \frac{-ac}{3m+7} \quad \cdots \quad \frac{-ac}{im+2i+1},   
\]
further also in the cases, in which $b$ is a term of this progression

\[
\frac{ac}{m+2}, \quad \frac{ac}{2(m+2)}, \quad \frac{ac}{3(m+2)} \quad \cdots \quad \frac{ac}{im+2i}.  
\]
But this continued fraction   exhibits the integral of the propounded differential equation with the boundary condition, that for $x=0$ it is $cxy=1$, if $m+2>0$; but if $m+2<0$, then the boundary condition is, that for $x=\infty$ it is $cxy=1$.

\paragraph*{§77}

Let us suppose that $b=0$ and $a=nc$ and that after the integration $x$ is put $=1$; from this equation

\[
ncx^mdx+cy^2dx+dy=0
\]
the following continued fraction will result, by which the value of $y$ will be defined in the case, in which one puts $x=1$,

\[
y=\cfrac{1}{c}+\cfrac{n}{\cfrac{-(m+3)}{c}+\cfrac{n}{\cfrac{2m+5}{c}+\cfrac{n}{\cfrac{-(2m+7)}{c}+\cfrac{n}{\cfrac{4m+9}{c}+\text{etc.}}}}};
\]
or just put $c=\frac{1}{\varkappa}$; from the equation

\[
nx^mdx+y^2dx+\varkappa dy=0
\]
the value of $y$ in the case, in which $x=1$, will be as follows

\[
y=\varkappa +\cfrac{n}{-(m \varkappa+3 \varkappa) +\cfrac{n}{2m \varkappa+5\varkappa +\cfrac{n}{-(3m \varkappa +7 \varkappa)+\text{etc.}}}}
\]
or 

\[
y=\varkappa -\cfrac{n}{m \varkappa+3 \varkappa -\cfrac{n}{2m \varkappa+5\varkappa -\cfrac{n}{3m \varkappa +7 \varkappa-\cfrac{n}{4m \varkappa +9 \varkappa -\text{etc.}}}}}
\]

\paragraph*{§78}

Therefore, if this continued fraction is propounded

\[
b+\cfrac{1}{b+\beta+\cfrac{1}{b+2 \beta+\cfrac{1}{b+3 \beta +\cfrac{1}{b+4 \beta +\text{etc.}}}}}
\]
it will be

\[
\varkappa=b, \quad n=-1, \quad (m+2)b=\beta
\]
or

\[
m=\frac{\beta}{b}-2.
\]
Therefore, the value of this continued fraction will be the solution for $y$ in the case, in which it is  $x=1$, of this differential equation

\[
x^{\frac{\beta -2b}{b}}dx=y^2dx+bdy
\]
having integrated it in such a way that for $x=0$ $xy$ becomes $=b$, because it is

\[
m+2>0,
\]
if $\frac{\beta}{b}$ is a positive number, of course.
\end{document}